\documentclass[12pt]{amsart}
\usepackage[T1]{fontenc}
\usepackage{times,mathptm}
\usepackage{amssymb,verbatim}
\usepackage[all]{xy}
\usepackage[dvips]{graphicx}
\usepackage{psfrag}
\usepackage{wrapfig}
\usepackage{epsfig,subfigure,color}
\usepackage{a4wide}

\newcommand{\cita}{\scshape}
\newcommand\Dq{q}
\newcommand{\R}{\mathbb{R}}

\newcommand{\Q}{\mathcal{Q}}
\newcommand{\HA}{\mathcal{H}}

\newcommand\D{d\hspace{-0.5pt}}
\newcommand\pA{\phi}
\newcommand{\dil}{\theta}

\theoremstyle{plain}
\newtheorem{Theorem}{Theorem}[section]
\newtheorem{Corollary}[Theorem]{Corollary}
\newtheorem{Proposition}[Theorem]{Proposition}
\newtheorem{Lemma}[Theorem]{Lemma}
\newtheorem{Remark}{Remark}[section]

\newtheorem{Definition}{Definition}[section]

\newtheorem*{NoNumberTheorem}{Theorem}
\newtheorem*{NoNumberProposition}{Proposition}

\newtheorem*{Convention}{Convention}

\begin{document}

\title[Pseudo-Anosov on hyperelliptic surfaces]
{Pseudo-{\cita Anosov} homeomorphisms on translation surfaces in hyperelliptic 
components have large entropy}

\author{Corentin Boissy, Erwan Lanneau}
\address{
LATP, case cour A,
Facult\'e de Saint J\'er\^ome
Avenue Escadrille Normandie-Niemen,
13397 Marseille cedex 20 }
\email{corentin.boissy@latp.univ-mrs.fr}

\address{
Centre de Physique Th\'eorique (CPT), UMR CNRS 6207 \newline
Universit\'e du Sud Toulon-Var and \newline
F\'ed\'eration de Recherches des Unit\'es de 
Math\'ematiques de Marseille \newline
Luminy, Case 907, F-13288 Marseille Cedex 9, France
}

\email{lanneau@cpt.univ-mrs.fr}


\begin{abstract}
We prove that the dilatation of any pseudo-{\cita Anosov} homeomorphism on a translation surface 
that belong to a hyperelliptic component is bounded from below uniformly by $\sqrt{2}$. This is in 
contrast to {\cita Penner}'s asymptotic. {\cita Penner} proved that the logarithm of the least dilatation of {\it any}
pseudo-{\cita Anosov} homeomorphism tends to zero at rate $1/g$ (as $g$ goes to infinity).

We also show that our uniform lower bound $\sqrt{2}$ is sharp. The proof 
uses the {\cita Rauzy-Veech} induction.
\end{abstract}

\date{May 22, 2010}
\maketitle
\setcounter{tocdepth}{1}
\section{Introduction}

Pseudo-{\cita Anosov} homeomorphisms play  an important role in {\cita
Teichm\"uller}  theory.  Calculating  problems  for their  dilatations
have  a long  history in  differential geometry.   The  unit cotangent
bundle  of the  moduli  space  of compact  genus  $g$ {\cita  Riemann}
surfaces $\mathcal M_g$ can be viewed as the moduli space of quadratic
differentials $\mathcal Q_g \rightarrow  \mathcal M_g$.  This space is
naturally  stratified  by   strata  of  quadratic  differentials  with singularities of
prescribed  multiplicities. The  {\cita  Teichm\"uller} geodesic  flow
$g_{t}$  acts naturally  on these  strata, and  closed loops  of length
$\log(\dil)>0$  for  this  flow  correspond to  conjugacy  classes  of
Pseudo-{\cita  Anosov}  homeomorphisms  with dilatation  $\dil>1$.  An
important  problem concerns  the asymptotic  behavior of  the smallest
dilatations. \medskip

The strata are not necessarily connected (see {\cita Kontsevich \&
Zorich}~\cite{KZ} and  {\cita Lanneau}~\cite{Lanneau:ccqd}). In  this
paper we use a discretization of the {\cita Teichm\"uller} geodesic
flow, i.e. the {\cita Rauzy-Veech} induction,  in order to tackle the
minimisation problem for {\it hyperelliptic components} $\mathcal
C^{hyp}$. This is  the first instance of asymptotic for components of
the moduli spaces. We shall prove

\begin{NoNumberTheorem}
Let $g\geq 1$. Let $\delta(\mathcal C^{hyp})$ be the least dilatation
of pseudo-{\cita Anosov} homeomorphisms affine on translation surfaces that belong to
$\mathcal C^{hyp} \subset \Q_{g}$. Then
$$ \sqrt{2} < \delta(\mathcal C^{hyp}) < \sqrt{2} +
\frac{4}{\sqrt{2}^{g}}.
$$
\end{NoNumberTheorem}

We will give a more precise statement in the following.

\subsection*{Mapping class group}

The  mapping  class group  $\textrm{Mod}(S)$  of  a closed  orientable
surface $S$ of genus $g \geq 1$ is defined to be the group of homotopy
classes   of  orientation  preserving   homeomorphisms  of   $S$.   An
irreducible  mapping class is  an isotopy  class of  homeomorphisms so
that no powers preserve a nontrivial subsurface of $S$.  By the {\cita
Thurston-Nielsen}    classification~\cite{Thurston1988},   irreducible
mapping classes are  either periodic (analogous to roots  of unity) or
are  of  a  type  called pseudo-{\cita  Anosov}~\cite{FLP}.   To  each
pseudo-{\cita  Anosov} mapping class  $[\pA] \in  \textrm{Mod}(S)$ one
can attach a  dilatation factor $\dil(\pA) > 1$. The logarithm  of $\dil(\pA)$
can be viewed as the minimal topological  entropy of any  element in the homotopy  class of
$\pA$ (uniquely realized by some element, $\pA$). \smallskip

{\cita Thurston} proved  that this number is an  algebraic integer and
even a  {\cita Perron}  number. $\dil(\pA) > 1$ is also the exponential  growth rate of 
lengths of curves under  iteration of $\pA$ (in any metric on $S$).   
These  numbers appear naturally  as the length spectrum  of the
moduli space of genus $g$ {\cita Riemann} surfaces. \medskip

It  is an  open question  to characterize  the set  of  dilatations of
pseudo-{\cita   Anosov}    homeomorphisms.    {\cita   Thurston}   has
conjectured  that pseudo-{\cita  Anosov} dilatations  (ignoring genus)
are precisely  the algebraic  units that are  {\cita Perron}  and also
larger than the {\cita Galois} conjugates of their inverses.

\subsection*{Minimisation problem}
Since    the   set    of    dilatations   (for    fixed   genus)    is
discrete~\cite{Arnoux:Yoccoz,Ivanov} (as  a subset of  $\R$) the least
dilatation $\delta_{g}$ is well defined. We know very little about the
values of the constants  $\delta_{g}$. The precise value of $\delta_2$
has been  recently calculated  (see {\cita Cho  \& Ham}~\cite{Cho:Ham}
and {\cita  Lanneau \& Thiffeault}~\cite{Lanneau:Thiffeault})  but the
values of $\delta_{g}$ for $g\geq 3$ are still unknown.  \medskip

Upper bounds are not hard to  derive from examples and there are a lot
of  results  in that  direction  (see  e.g.  {\cita Penner,  McMullen,
Hironaka, Kin, Minakawa}). So far  the best general upper bound is the
one   given   by    {\cita   Hironaka   \&   Kin}~\cite{Hironaka2006}:
$\log(\delta_g)  \leq  \frac{\log(2+\sqrt{3})}{g}$. See also~\cite{Hironaka2009,K2010} for 
recent results in that direction. But  again  very
little    is    known     about    {\it    lower}    bounds.    {\cita
Penner}~\cite{Penner1991}    proved    that    $\log(\delta_g)    \geq
\frac{\log(2)}{12g-12}$,  using  general   properties  of  the  {\cita
Perron-Frobenius}  matrices.   There  is  also  a   result  of  {\cita
Tsai}~\cite{Tsai2009}  for pseudo-{\cita  Anosov}  on {\it  punctured}
surfaces. \medskip

In  general  lower  bounds  are  much  subtle  to  obtain  than  upper
bounds.  In  contrast  to  our  understanding  of  the  asymptotic  of
$\log(\delta_{g})$, we still  do not know the answer  to the following
question, posed by {\cita McMullen}~\cite[Section~10]{Mc}:
$$  \textrm{Does }  \lim\limits_{g \to  +\infty} g\log(\delta_{g})
\textrm{ exist? What is its value?}
$$

\subsection*{Subgroups of the modular group and strata}
In his book~\cite{Farb}, {\cita Farb} proposed two natural refinements
of the minimisation problem. \medskip

The  first one  is related  to subgroups  of the  modular  group. More
precisely, let us fix a subgroup $H \subseteq \textrm{Mod}(S)$ and let
us consider the least  dilatation $\delta(H)$ of pseudo-{\cita Anosov}
classes $[\pA]\in H$. \medskip

In order to state the second  problem, let us recall the definition of
strata here (see Section~\ref{sec:background} for precise statements).
A pseudo-{\cita  Anosov} homeomorphism on a compact  surface defines a
pair of transverse measured foliations. For simplicity, we will assume
that  the  foliations  are  orientable,  but the  techniques  we  will
develop   holds   for   non-orientable  measured   foliations   (see
Section~\ref{sec:quadratic}).  These  data are  equivalent to the  given of  a pair
$(X,\omega)$ (called a  translation surface for a reason  that will be
clear in  the next) where $X$ is  a Riemann surface and  $\omega$ is a
holomorphic     $1-$form     with     zeroes     of     multiplicities
$(k_1,\dots,k_n)$. We will also say singularities of multiplicities $(k_1,\dots,k_n)$ of the 
translation surface.  The   {\cita    Gauss-Bonnet}    formula   reads
$\sum_{i=1}^n  k_i  =  2g-2$.  A   stratum  consists  of  the  set  of
translation surfaces  with prescribed singularities.   {\cita Masur \&
Smillie}~\cite{Masur:Smillie} proved that, for every singularity data,
there exits some pseudo-{\cita  Anosov} homeomorphism that realize it.
Hence it  is natural to  consider the following  problem (minimisation
problem in a stratum)
$$ \delta(k_1,\dots,k_n) := \inf \left\{
\log(\pA):
\begin{array}{l}
\textrm{pseudo-{\cita  Anosov} $\pA$  whose  corresponding} \\ \textrm{translation surface has singularity data } (k_1,\dots,k_n)
\end{array}
\ \right\}. $$ 
Similarly, for a pair of nonoriented measured foliation, one gets a pair $(X,q)$, where $X$ is a Riemann surface and $q$ is a holomorphic quadratic differential with zeroes of order $(k_1,\ldots ,k_n)$.
Obviously one has
$$ \delta_g = \min \{ \delta(k_1,\dots,k_n) \}, $$ 
where the min is taken over all possible singularity data (with no assumptions on the 
orientability of the foliations). \medskip

It turns out that strata  (translation surfaces with prescribed singularity data
$(k_1,\dots,k_n)$)  are not
connected in general~\cite{KZ,Lanneau:ccqd}.  We can thus consider the
natural following  refinement of {\cita  Farb}'s question~: minimizing
problem for connected components. \\
The components are distinguished by  two invariants the {\it parity of
the spin  structure} and the  {\it hyperelliptic components}.  The  latter are defined as
follows.

Recall  that  a  {\cita  Riemann}   surface $X$  of  genus  $g\geq  2$  is
hyperelliptic if  there exists  an holomorphic involution $\tau$  with $2g+2$
fixed points.
\begin{Convention}
In all  of this paper,  we will use  the following convention:  a flat
surface  $(X,\omega)$ is  hyperelliptic if  the underlying  {\cita Riemann}
surface  is  hyperelliptic  and   $\tau^*\omega=-\omega$  (or
equivalently      is       an      affine      homeomorphism,      see
Section~\ref{sec:homeo}). The  fixed points of $\tau$ will be
usually called {\cita Weierstrass} points.
\end{Convention}

Hyperelliptic connected components $\mathcal{C}^{hyp}$ are defined in the following way

\begin{equation*}
\label{eq:hyp}
\mathcal C^{hyp}:= \left\{ (X,\omega):
\begin{array}{l}
X \textrm{  is hyperelliptic; }  \omega \textrm{ is a  holomorphic one
form   with  one  single}   \\  \textrm{zero,   or  two   zeroes;  the
hyperelliptic  involution preserves }  \omega \\  \textrm{and permutes
the zeroes, if there are two}
\end{array}
\right\}.
\end{equation*}

{\cita Leininger}~\cite{Leininger2004} and then {\cita Farb, Leininger
\& Margalit}~\cite{Farb2008}  tackle the minimisation  problem for the
subgroups of $Mod(S)$ given by the {\cita Thurston}'s construction and
for the {\cita Torelli} group, respectively. They provide evidence for
the   principle   that    algebraic   complexity   implies   dynamical
complexity. \medskip

In  this paper  we will  investigate the  second problem  and  prove a
similar    theorem.     This   answers    a    question   of    {\cita
Farb}~\cite[Problem 7.5]{Farb}. We shall prove

\begin{Theorem}
\label{theo:main:general}
Let $g\geq 1$. Let $\pA$ be a pseudo-{\cita Anosov} homeomorphism affine on a
genus $g$ translation surface $(X,\omega) \in \mathcal C^{hyp}$.  Then
$$ \dil(\pA) > \sqrt{2}.
$$
\end{Theorem}

Recall that {\cita Penner}~\cite{Penner1991} proved that, as the genus
increases,  there   are  pseudo-{\cita  Anosov}   homeomorphisms  with
dilatations arbitrarily close to $1$. \medskip

\noindent This uniform lower bound  is sharp. Indeed one can construct
(see Appendix~\ref{appendix:examples})  a sequence $(\pA_g)_{g\geq 1}$
on  $(X,\omega)\in \mathcal  C^{hyp}$ with  dilatations  converging to
$\sqrt{2}$. More precisely, we will show

\begin{Theorem}
\label{theo:main:sharp}
Let $g\geq 1$. Let  $\delta(\mathcal C^{hyp})$ be the least dilatation
of  pseudo-{\cita Anosov}  homeomorphisms affine on  a genus  $g$ translation
surface $(X,\omega) \in \mathcal C^{hyp}$. Then
$$    \sqrt{2}    <   \delta(\mathcal    C^{hyp})    <   \sqrt{2}    +
\frac{4}{\sqrt{2}^{g}}.
$$
\end{Theorem}

\subsection*{Non hyperelliptic components}
The reader may wonder why we  impose the restrictions on the action of
the  hyperelliptic  involution on  the  zeroes  in  the definition  of
$\mathcal C^{hyp}$. It  turns out that if we don't impose this restriction, then the translation surface is hyperelliptic, but not in the hyperelliptic connected component, and therefore, the property of being hyperelliptic is not preserved by any small deformation inside the ambiant stratum. We will show 
 that the asymptotic  behavior may be
very different if we consider different connected components. \medskip

As       for      example       one      can       construct      (see
Appendix~\ref{appendix:examples:2}) a sequence $(\varphi_g)_{g\geq 3}$
(with  $g$   odd)  of   pseudo-{\cita  Anosov}  homeomorphisms   on  a
hyperelliptic  translation surface  $(X_{g},\omega_{g})$ of  genus $g$
having two zeroes of degree $g-1,g-1$, and such that the dilatation of
$\varphi_g$ is the {\cita Perron} root of the polynomial
$$ X^{2g} - X^{2g-1} - 4X^g - X + 1.
$$ In particular
$$ \lim\limits_{k \to +\infty} \dil(\varphi_{2k+1}) = 1.
$$ Of course  in that case the hyperelliptic  involution fixes the two
zeroes. \medskip

In  the  above  case,  the ({\it  non-hyperelliptic})  components  are
distinguished  by a parity  of the  spin structure  (see~\cite{KZ} and
Section~\ref{sec:components}).  Let $\mathcal  C^{odd}$  be the  ({\it
non-hyperelliptic}) odd component. Then we shall prove

\begin{Theorem}
\label{theo:main:odd}
Let $g\geq 1$, be an odd integer. Let  $\delta(\mathcal C^{odd})$ be the least dilatation
of  pseudo-{\cita Anosov}  homeomorphisms affine on  a genus  $g$ translation
surface $(X,\omega) \in \mathcal C^{odd}$. Then
$$ 1 < \delta(\mathcal C^{odd}) < 1 + \frac1{g}.
$$
\end{Theorem}

\subsection*{Quadratic differentials}

A quadratic  differential $(X,q)$ is  strictly quadratic if it  is not
the  square  of  any  holomorphic $1$-form.  $q$ determines foliations with singularities ($p$-prongs). 
There exists a canonical branched double
cover  $Y \rightarrow X$  on which  $q$ becomes  a square.  The branch
points correspond to  the simple poles ($1$-prongs) and zeros of  odd order of $q$.
We    will   prove    Theorem~\ref{theo:main:general}    for   Abelian
differentials or equivalently translation surfaces. \medskip

\noindent  In Section~\ref{sec:quadratic} we  derive from  our theorem
some results for strict  quadratic differentials, namely we will prove
Theorem~\ref{theo:quadratic}, page~\pageref{theo:quadratic}. \medskip

\noindent As mentioned, one has a more general result, with no assumptions whether the 
foliations are orientable or not.

\begin{Theorem}
\label{theo:main:general:2}
Let $g\geq 1$. Let $\pA$ be a pseudo-{\cita Anosov} homeomorphism on  a genus  $g$ 
hyperelliptic surface. Let us assume that the foliations have some 1-prongs and only one or two 
p-prongs with $p\geq 2$. Assume also that the hyperelliptic involution is affine (for the flat metric). If the involution permutes the 
singularities (if there are two) then
$$    
\dil(\pA) > \sqrt{2}.
$$
\end{Theorem}

\subsection*{{\cita Rauzy-Veech} induction}

If $\mathcal C$ is a connected  component of some strata then there is
a finite ramified  covering $\widehat{\mathcal C} \rightarrow \mathcal
C$  consisting of  marking  a zero  and  a separatrix  of surfaces  in
$\mathcal C$.  The {\cita  Rauzy-Veech} induction provides  a discrete
representation (symbolic coding) of  the {\cita Teichm\"uller} flow on
$\widehat{\mathcal              C}$~\cite{Veech1982}              (see
Section~\ref{appendix:construction:pA}   for   precise   definitions).
Periodic  orbits in  $\mathcal C$  are taken  to conjugacy  classes of
pseudo-{\cita Anosov} on $(X,\omega) \in \mathcal C$. \medskip

\noindent  As   shown  by  {\cita  Veech},  to   each  periodic  orbit
$\gamma\subset \widehat{\mathcal  C}$ there corresponds  a closed loop
in   some   graph  (called   a   {\cita   Rauzy}  diagram)   $\mathcal
D_{r}(\mathcal  C)$ and a  renormalization matrix  $V(\gamma)\in SL(h,
{\mathbb  Z})$ ($h=2g+n-1$).   This  matrix  corresponds  to  the   action  of  the
corresponding pseudo-{\cita Anosov} homeomorphism in relative homology
of the  underlying surface  with respect to  the singularities  of the
Abelian differential. Hence the  spectral radius of $V(\gamma)$ is the
dilatation  of  the  pseudo-{\cita  Anosov}  on  $(X,\omega)  \in  \mathcal
C$.\medskip

\noindent A crucial point in this paper is a carefully analysis of the
geometry of  these {\cita Rauzy} diagrams.  This  approach was already
used by  {\cita Avila  \& Viana}~\cite{Avila:Viana} for  the dynamical
properties of the {\cita Teichm\"uller} geodesic flow.

\subsection*{Outline of a proof of our main result}

The {\cita  Rauzy-Veech} induction allows one  to relate pseudo-{\cita
Anosov} homeomorphisms and closed  loops in {\cita Rauzy} diagrams. We
conclude    by    sketching    its     use    in    the    proof    of
Theorem~\ref{theo:main:general}.

\begin{enumerate}

\item[1.] Let  $\pA$ be  a pseudo-{\cita Anosov}  homeomorphism affine
with  respect  to  a  translation  surface  $(X,\omega)  \in  \mathcal
C^{hyp}$.  We  prove  that   $\pA$  commutes  with  the  hyperelliptic
involution  (Proposition~\ref{prop:commutes}).  Thus  $\pA$  induces a
pseudo-{\cita Anosov} homeomorphism on the sphere.

\item[2.] Using the {\cita Brouwer}  fixed point theorem, we show that
$\pA^2$ fixes a separatrix of the horizontal measured foliation on $X$
(Proposition~\ref{theo:reduc1}).

\item[3.] Hence $\pA^{2}$ is also affine with respect to a surface $(X,\omega) \in
\widehat{\mathcal C^{hyp}}$ (with  a marked separatrix).  Thus $\pA^2$
is obtained by taking a (irreducible) closed loop in the {\cita Rauzy}
diagram  $\mathcal  D_{r}^{hyp}$  corresponding to  $\widehat{\mathcal
C^{hyp}}$.

\item[4.]   The   transition  matrices   of   the  diagram   $\mathcal
D_{r}^{hyp}$  are  hard to  describe.   We  use {\cita Marmi-Moussa-Yoccoz}'s
representation   which   furnished   a  covering   $\mathcal   D^{hyp}
\rightarrow \mathcal  D_{r}^{hyp}$. We then obtain  a simple criterion
for a path $\gamma$ (not necessarily  closed) in $\mathcal D^{hyp}$ to
have     a transition matrix $V(\gamma)$ with    spectral       radius        greater        than       $2$
(Proposition~\ref{prop:win:los}).

\item[5.] Analyzing carefully the combinatorics of these {\cita Rauzy}
diagrams              $\mathcal              D^{hyp}$, we  show that $\dil(\pA^2) \geq
2$ (see
Section~\ref{sec:description:marked}). Since $\dil(\pA^2) = \dil(\pA)^2$ one gets the desired result.

\item[6.]  In  Appendix~\ref{appendix:examples}   we  prove  that  our
uniform  bound is  sharp  by exhibiting  a  sequence of  pseudo-{\cita
Anosov} homeomorphisms. In Appendix~\ref{appendix:examples:2} we show that the action
of the hyperelliptic involution on the zeroes is important.

\end{enumerate}

\subsection*{Further results}

{\cita   Fehrenbach  \&   Los}~\cite{Los1997}  proved   the  following
inequality  for pseudo-{\cita  Anosov} homeomorphisms  $\pA$  having a
periodic   orbit  of   length  $n\geq   3$:   $\textrm{log}(\pA)  \geq
\frac1{n}\textrm{log}(1  + \sqrt{2})$.  However  this does  not easily
imply a  uniform lower bound. Indeed for  each hyperelliptic connected
component,  there exists  a pseudo-{\cita  Anosov}  homeomorphism such
that the minimal length (greater than  two) of a periodic orbit is the
number  of {\cita  Weierstrass}  points i.e.  $2g+2$  points (see  the
appendix). In particular, Theorem~\ref{theo:main:general} does not follow from these
estimates.

 \medskip

See   also~\cite{Los2009,Lanneau:fixed}   for   further   results   on
pseudo-{\cita  Anosov}  that do  not  come  from  the {\cita  Veech}'s
construction.

\subsection*{Acknowledgements}

We thank Vincent {\cita Delecroix} for remarks and comments on this paper and for his
very useful sage's library ``IET''~\cite{sage}.

\section{Background}
\label{sec:background}

We   review   basic   notions   and  results   concerning   Abelian
differentials, translation surfaces, pseudo-{\cita Anosov} homeomorphisms and moduli spaces. For general
references                                                          see
say~\cite{MT,Veech1982,Rauzy,Masur:Smillie,Marmi:Moussa:Yoccoz}.

\subsection{Flat surfaces and pseudo-{\cita Anosov} homeomorphisms}
\label{sec:homeo}

A {\it  flat surface}  $S$ is a  (real, compact, connected)  genus $g$
surface  equipped  with  a  flat  atlas i.e.   a  triple  $(S,\mathcal
U,\Sigma)$  such  that $\Sigma$  is  a  finite  subset of  $S$  (whose
elements   are  called   {\em  singularities})   and  $\mathcal   U  =
\{(U_i,z_i)\}$  is an atlas  of $S  \setminus \Sigma$  with transition
maps $z  \mapsto \pm\ z +  \textrm{ constant}$.  We  will require that
for each $s\in \Sigma$, there is  a neighborhood of $s$ isometric to a
Euclidean  cone.   Therefore we  get  a  {\it quadratic  differential}
defined  locally  in the  coordinates  $z_i$  by  the formula  $\Dq=\D
z_i^2$. This form extends to  the points of $\Sigma$ to zeroes, simple
poles or marked  points (we will usually call the  zeroes and poles \emph{singular points} or  simply \emph{singularities}). \medskip

If there  exists a  sub-atlas such that  all transition  functions are
translations  then  the quadratic  differential  $\Dq$  is the  global
square of  an Abelian differential $\omega \in  H^1(X,\mathbb C)$.  We
will then say that $(X,\omega)$ is a translation surface.  \medskip

A homeomorphism $f : X \rightarrow X$ is an {\it affine homeomorphism}
if  $f$ restricts to  a diffeomeomorphism  of $X\setminus  \Sigma$ of
constant derivative.  It is equivalent to say that $f$ restricts to an
isomorphism of $X\setminus \Sigma$ which preserves the induced affine structure given
by $q$. \medskip

There    is    a    standard    classification    of    elements    of
$\textrm{SL}_2(\mathbb R)$  into three types:  elliptic, parabolic and
hyperbolic. This  induces a classification  of affine diffeomorphisms.
An affine  diffeomorphism is parabolic, or  elliptic, or pseudo-{\cita
Anosov},    respectively,   if    $|\textrm{trace}(D    f)|   =    2$,
$|\textrm{trace}(D  f)|   <  2$,  or  $|\textrm{trace}(D   f)|  >  2$,
respectively.  If  $\pA$ is pseudo-{\cita Anosov},  in the coordinates
of the  stable and unstable  measured foliations determined  by $\pA$,
one has $D\pA= \left( \begin{smallmatrix} \dil^{-1} & 0 \\ 0 & \dil
\end{smallmatrix} \right)$  where $|\dil| > 1$. The number $|\dil|$ is
called  the {\it  dilatation} of  $\pA$.  From now  all flat  surfaces
considered    will     be    translation    surfaces,     except    in
Section~\ref{sec:quadratic} and Appendix~\ref{appendix:examples:2}.

\begin{Remark}
Constructing  parabolic elements  is  easy (see~\cite{Veech1989})  but
constructing  hyperbolic  elements  is  more  subtle.  A  first
construction  is to consider  the product  of parabolic  elements (see
e.g.    the   paper    of    {\cita   Fathi}~\cite{Fathi1987}).     In
Section~\ref{appendix:construction:pA}   we  recall  a   very  general
construction due to {\cita Veech}.
\end{Remark}

For $g  \geq 1$, we define  the moduli space  of Abelian differentials
$\HA_g$ as the moduli space of pairs $(X,\omega)$ where $X$ is a genus
$g$  (compact,  connected)  {\cita  Riemann}  surface  and  $\omega\in
\Omega(X)$ a  non-zero holomorphic $1-$form defined on  $X$.  The term
moduli  space  means that  we  identify  the  points $(X,\omega)$  and
$(X',\omega')$   if  there   exists  an   analytic   isomorphism  $f:X
\rightarrow X'$ such that $f^\ast \omega'=\omega$.  The  group $\textrm{SL}_2(\mathbb  R)$  naturally acts  on the  moduli space  of flat  surfaces by  post composition  on the  charts.  \medskip

One can also see a translation surface obtained as a polygon (or a union of polygons) whose sides come by pairs, and for each pairs, the corresponding segments are parallel and of the same lengths. These parallel sides are glued together by translation and we assume that this identification preserves the natural orientation of the polygons. In this context, two translation surfaces are identified in the moduli space of Abelian differentials if and only if the corresponding polygons can be obtained from each other by ``cutting'' and ``gluing'' and preserving the identifications (i.e.  the two surfaces represent the same point in the moduli space). Also, the $SL_2(\mathbb{R})$ action in this representation is just the natural linear action on the polygons. \medskip 

{\cita Veech} showed that an affine homeomorphism with a derivative map which is  not  the  identity is  not  isotopic  to  the identity.  Hence  a
homeomorphism  $f$ is  an  affine homeomorphism  on  the flat  surface
$(X,\omega)$, with  derivative map  $Df=A$, if and  only if the  matrix $A$
stabilizes the surface  $(X,\omega)$. That is $(X,\omega)$ can  be obtained from $A \cdot (X,\omega)$  by ``cutting'' and ``gluing''.

\subsection{Connected components of the strata}
\label{sec:components}

The  moduli  space  of  Abelian  differentials is  stratified  by  the
combinatorics of  the zeroes;  We will denote  by $\HA(k_1,\dots,k_r)$
the stratum  of $\HA_g$ consisting of (classes  of) pairs $(X,\omega)$
such  that   $\omega$  possesses  exactly  $r$  zeroes   on  $X$  with
multiplicities $(k_1,\dots,k_r)$.

It is a well known part of the {\cita Teichm\"uller} theory that these
spaces  are   ({\cita  Hausdorff})  complex  analytic,   and  in  fact
algebraic, spaces.

These strata are non-connected in general but each stratum has at most
three   connected    components   (see~\cite{KZ}   for    a   complete
classification). In particular for $g\geq 4$ the stratum with a single
zero,  $\HA(2g-2)$,  has   three  connected  components.  The  stratum
$\HA(g-1,g-1)$ has two or three connected components depending whether
$g$ is even or odd, respectively.

\subsubsection{Hyperelliptic component}
\label{sec:hyp}
This component  contains precisely pairs  $(X,\omega)$ where $X$  is a
hyperelliptic surface  and $\omega$ a one-form whose  zeroes (if there
are  two)  are  interchanged   by  the  hyperelliptic  involution.  An
equivalent  formulation is  to require  that there  exists  a ramified
double cover $\pi  : X \rightarrow \mathbb P^1$ over  the sphere and a
quadratic differential $q$  on $\mathbb P^1$ having only  one zero and
simples poles  such that $\omega^2  = {\pi^{*} q}$. We  will denote
these components by $\HA^{hyp}(2g-2)$ and $\HA^{hyp}(g-1,g-1)$.

\subsubsection{Other components}
The other ({\it non-hyperelliptic})  components are distinguished by a
parity  of the  spin structure.   There are  two ways  to  compute the
parity of the  spin structure of a translation  surface $X$. The first
way  is  to  use  the  {\cita  Arf}  formula  on  a  symplectic  basis
(see~\cite{KZ}). The  second possibility applies  if $X$ comes  from a
quadratic differential, i.e.\ if $X$ possesses an involution such that
the         quotient        produces         a        half-translation
surface~\cite{Lanneau:spin}.     We     will     apply     this     in
Appendix~\ref{appendix:examples:2}.

\subsection{Application of Theorem~\ref{theo:main:general} to quadratic differentials}
\label{sec:quadratic}

In  this section,  we extend  Theorem~\ref{theo:main:general}  to some
other strata in the moduli space of quadratic differentials. This part
is independent  from the rest  of the paper  and can be skipped  for a
first reading. However, Proposition~\ref{prop:commutes} will be needed
later. \medskip

As for Abelian differentials, strata  of the moduli space of quadratic
differentials  are not  connected in  general (see~\cite{Lanneau:ccqd}
for    a    complete    classification).    We   can    deduce    from
Theorem~\ref{theo:main:general} results on hyperelliptic components in
the quadratic case. \medskip

We denote  by $\Q(k_1,\dots,k_n)$ strata  of the moduli space  of half
translation surfaces where the vector $(k_1,\dots,k_n)$ agrees with the {\cita Gauss-Bonnet} formula 
$\sum_{i=1}^n k_i  = 4g-4$. For  $g\geq 2$,
let   us   consider   the   two   strata   $\Q(-1,-1,2g-3,2g-3)$   and
$\Q(-1,-1,4g-2)$.  Their  hyperelliptic component  can  be defined  as
follows (see~\cite{Lanneau:hyp}):
$$ \mathcal C_1^{hyp}:= \left\{ (X,q):
\begin{array}{l}
X \textrm{ is  hyperelliptic, } (X,q)\in \Q(-1,-1,2g-3,2g-3), \textrm{
and} \\  \textrm{the involution  permutes the two  zeroes and  the two
poles}
\end{array}
\right\}.
$$
$$ \mathcal C_2^{hyp}:= \left\{ (X,q):
\begin{array}{l}
X \textrm{ is hyperelliptic,  } (X,q)\in \Q(-1,-1,4g-2), \textrm{ and}
\\ \textrm{the involution permutes the two poles}
\end{array}
\right\}.
$$

\begin{Theorem}
\label{theo:quadratic}
Let  $g\geq 2$.  Let  $\mathcal C$  be  one of  the two  hyperelliptic
components  defined  above.   Let  $\pA$ be  a  pseudo-{\cita  Anosov}
homeomorphism  affine on  a   half  translation  surface  $(X,q)\in  \mathcal
C$. Then  the dilatation $\dil(\pA)$ satisfies  $\dil(\pA) >
\sqrt{2}$.
\end{Theorem}

We will use the following result.

\begin{Proposition}
\label{prop:commutes}
Let $f$  be an affine  homeomorphism on hyperelliptic flat  surface of
genus greater than or equal to two.
The map $f$ commutes with the hyperelliptic involution.
\end{Proposition}

\begin{proof}[Proof of Proposition~\ref{prop:commutes}]
Let  $\tau$ be the  hyperelliptic involution  on $X$.   The commutator
$[f,\tau]  =  f  \circ   \tau\circ  f^{-1}\circ  \tau$  is  an  affine
homeomorphism  and  its  derivative  map  is the  identity  matrix  in
$\textrm{SL}_2(\R)$. So $[f,\tau]$ is a translation. 

Since the genus is greater than or equal to two, and the surface is hyperelliptic, then the set of  {\cita Weierstrass} points is nonempty.
We  show  that  $f$  preserves  set-wise  the  set  of  {\cita
Weierstrass} points  (i.e. the set $W =  \textrm{Fix}(\tau)$). The map
$f^{-1}\tau f$ is a conformal  automorphism of the complex surface, it
is  easy  to see  that  it is  a  hyperelliptic  involution. Hence  it
preserves        point-wise       the        {\cita       Weierstrass}
points~\cite{Farka}. \medskip Let $p \in W$. Then, $f^{-1}\tau f(p)=p$
so  $\tau f(p)=f(p)$. Hence  $f(p)$ is  a fixed  point of  $\tau$ thus
$f(p)  \in  W$.  Since  $f$  preserves set-wise  $W$  the  translation
$[f,\tau]$  fixes {\it  point-wise} $W$.  It is  easy to  see  that it
preserves also  the separatrices  issued from the  regular {\cita Weierstrass}
points, that necessarily exists.  Thus   the  commutator  $[f,\tau]$  is   the  identity.   The
proposition is proven.
\end{proof}

\begin{proof}[Proof of Theorem~\ref{theo:quadratic}]
Let $(X,q)\in  \Q(-1,-1,2g-3,2g-3)$ be  a half translation  surface in
the  hyperelliptic component.  Let  $\pA$ be  a pseudo-{\cita  Anosov}
homeomorphism on $X$ and $\tau$ the hyperelliptic involution. \medskip

Passing to  the quotient we  get a meromorphic  quadratic differential
$q'$ on the projective line $\mathbb P^1$ with $2g+1$ simple poles and
a single zero of degree $2g-3$.  Taking the standard orientating cover
$\pi:Y \rightarrow \mathbb P^1$ over $\mathbb P^1$ having ramification
points precisely  over odd degree singularities (namely  the poles and
the zero) we obtain a translation surface $(Y,\omega)$ where $\omega^2 =
{\pi^{*}  q'}$.  By  construction (see  Subsection~\ref{sec:hyp})
$(Y,\omega)$     belongs     to     the    hyperelliptic     component
$\HA^{hyp}(2g-2)$.   Now   by  Proposition~\ref{prop:commutes}   $\pA$
commutes  with  $\tau$ on  $Y$.  Thus  $\pA$  induces a  pseudo-{\cita
Anosov}  homeomorphism  $\varphi$  on  $\mathbb  P^1$  with  the  same
dilatation. Since  $\varphi$ preserves the set  of ramification points
of   $\pi$,   $\varphi$  lifts   to   a   new  pseudo-{\cita   Anosov}
homeomorphism, say $\hat \varphi$, on $Y$.
$$  \xymatrix{  X  \ar[r]^{\pA}  \ar[d]   &  X  \ar[d]  &  Y  \ar[lld]
\ar[r]^{\hat \varphi}  & Y  \ar[lld]\\ \mathbb P^1  \ar[r]^{\varphi} &
\mathbb P^1 }
$$  Now   $\dil(\pA)  =   \dil(\varphi)  =  \dil(\hat   \varphi)$.  By
Theorem~\ref{theo:main:general}    we    get    that   $\dil(\pA)    >
\sqrt{2}$. \medskip

The second case is similar and left to the reader.
\end{proof}

\noindent We end with the following (see the introduction): \medskip

\noindent {\bf Theorem~\ref{theo:main:general:2}.} 
{\it Let $g\geq 1$. Let $\pA$ be a pseudo-{\cita Anosov} homeomorphism on  a genus  $g$ 
hyperelliptic surface. Let us assume that the foliations have some 1-prongs and only one or two 
p-prongs with $p\geq 2$. Assume also that the hyperelliptic involution is affine (for the flat metric). If the involution permutes the 
singularities (if there are two) then
$$    
\dil(\pA) > \sqrt{2}.
$$
}

\begin{proof}[Proof of Theorem~ \ref{theo:main:general:2}]
Let $\pi : (X,\omega) \rightarrow (\mathbb P^{1},q)$ be the covering. By construction $\pA$ induces a 
pseudo-{\cita Anosov} homeomorphism on the sphere and lift to a pseudo-{\cita Anosov} homeomorphism 
on a translation surface in some hyperelliptic component, with the same dilatation. Then 
Theorem~\ref{theo:main:general} applies.
\end{proof}

\section{{\cita Rauzy-Veech} induction and pseudo-{\cita Anosov} homeomorphisms}
\label{appendix:construction:pA}

In  this section  we recall  the basic  construction  of pseudo-{\cita
Anosov}  homeomorphisms using the  {\cita Rauzy-Veech}  induction (for
details         see~\cite{Veech1982},         \S        8,         and
\cite{Rauzy,Marmi:Moussa:Yoccoz}).  We  first review the  link between
interval exchange maps and translation surfaces.

\subsection{Interval exchange transformations}
Let $I  \subset \mathbb  R$ be an  open interval  and let us  choose a
finite partition  of $I$  into $d\geq 2$  open subintervals  $\{I_j, \
j=1,\dots,d \}$.  An interval  exchange transformation is a one-to-one
map  $T$  from  $I$  to  itself that  permutes,  by  translation,  the
subintervals $I_j$. It is easy to see that $T$ is precisely determined
by the  following data:  a permutation that
encodes how  the intervals are  exchanged, and a vector  with positive
entries that encodes the lengths of the intervals. \medskip

We  will  use  the  description  given  by  {\cita  Marmi,  Moussa  \&
Yoccoz}~\cite{Marmi:Moussa:Yoccoz}. This will simplify the description
of  the induction  which will  be  very useful  for the  proof of  our
result.

We will  attribute a name to  each interval $I_{j}$. In  this case, we
will speak of \emph{labeled} interval  exchange maps.  One gets a pair
of  one-to-one   maps  $(\pi_t,\pi_b)$  (t  for  ``top''   and  b  for
``bottom'') from  a finite alphabet  $\mathcal{A}$ to $\{1,\ldots,d\}$
in  the following  way. In  the partition  of $I$  into  intervals, we
denote the  $k^{\textrm{th}}$ interval, when counted from  the left to
the right,  by $I_{\pi_t^{-1}(k)}$. Once the  intervals are exchanged,
the  interval  number  $k$  is  $I_{\pi_b^{-1}(k)}$.  Then  with  this
convention,  the permutation  encoding  the map  $T$  is $\pi_b  \circ
\pi_t^{-1}$. We  will denote the length  of the intervals  by a vector
$\lambda=(\lambda_\alpha)_{\alpha\in \mathcal{A}}$.

\begin{Definition}
We  will call  the  pair $(\pi_t,\pi_b)$  a  \emph{labeled} permutation,  and
$\pi_b   \circ    \pi_t^{-1}$   a   permutation    (or   \emph{reduced}
permutation). If it  is clear from the context, then  we will just use
the  term  permutation. We  will  also  usually  write a  reduced
permutation  as a  labeled one  with  $\mathcal{A}=\{1,\ldots,d\}$ and
$\pi_t=Id$.
\end{Definition}

One usually  represents labeled permutations $\pi=(\pi_t,\pi_b)$ by
a table:
\begin{eqnarray*}
\pi=
\left(\begin{array}{ccccc}\pi_t^{-1}(1)&\pi_t^{-1}(2)&\ldots&\pi_t^{-1}(d)
\\ \pi_b^{-1}(1)&\pi_b^{-1}(2)&\ldots&\pi_b^{-1}(d)
\end{array}\right).
\end{eqnarray*}

\subsection{Suspension data}
\label{sec:suspension}

The  next  construction  provides  a link  between  interval  exchange
transformations  and  translation surfaces.   A  suspension datum  for
$T=(\pi,\lambda)$      is       a      collection      of      vectors
$\{\tau_\alpha\}_{\alpha\in \mathcal{A}}$ such that
\begin{itemize}
\item  $\forall  1  \leq   k  \leq  d-1,\  \sum_{\pi_t(\alpha)\leq  k}
\tau_\alpha>0$,
\item  $\forall  1  \leq   k  \leq  d-1,\  \sum_{\pi_b(\alpha)\leq  k}
\tau_\alpha<0$.
\end{itemize}

We  will  often  use  the notation  $\zeta=(\lambda,\tau)$.   To  each
suspension  datum  $\tau$,  we  can associate  a  translation  surface
$(X,\omega)=X(\pi,\zeta)$ in the following way.

\begin{wrapfigure}[6]{r}{0pt}
\noindent
\psfrag{1}[][]{\small   $\zeta_1$}   \psfrag{2}[][]{\small  $\zeta_2$}
\psfrag{3}[][]{\small   $\zeta_3$}   \psfrag{4}[][]{\small  $\zeta_4$}
\includegraphics[width=6cm]{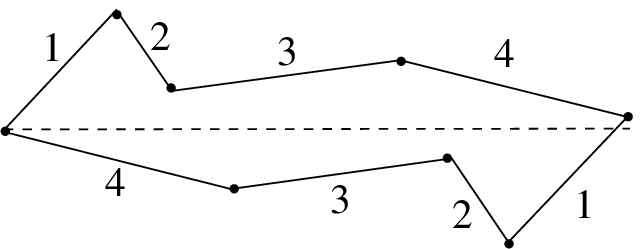}
\end{wrapfigure}
Consider the broken line  $L_t$ on $\mathbb{C}=\mathbb R^2$ defined by
concatenation of  the vectors $\zeta_{\pi_t^{-1}(j)}$  (in this order)
for  $j=1,\dots,d$ with starting  point at  the origin.  Similarly, we
consider the broken line $L_b$ defined by concatenation of the vectors
$\zeta_{\pi_b^{-1}(j)}$  (in   this  order)  for   $j=1,\dots,d$  with
starting point  at the origin.  If  the lines $L_t$ and  $L_b$ have no
intersections other than the endpoints, we can construct a translation
surface $X$ by identifying each  side $\zeta_j$ on $L_t$ with the side
$\zeta_j$  on $L_b$  by  a  translation. The  resulting  surface is  a
translation  surface endowed  with the  form $\D  z^2$. Note  that the
lines $L_t$ and $L_b$ might  have some other intersection points.  But
in this case, one can still  define a translation surface by using the
\emph{zippered   rectangle  construction},   due   to  {\cita   Veech}
(\cite{Veech1982}).

Let  $I  \subset  X$  be  the  horizontal interval  defined  by  $I  =
(0,\sum_{\alpha} \lambda_{\alpha})  \times \{0\}$.  Then  the interval
exchange transformation $T$ is precisely  the one defined by the first
return map of the vertical flow on $X$  to $I$.

\subsection{{\cita Rauzy-Veech} induction}
The {\cita Rauzy-Veech} induction $\mathcal R(T)$ of $T$ is defined as
the  first return  map of  $T$  to a  certain subinterval  $J$ of  $I$
(see~\cite{Rauzy,Marmi:Moussa:Yoccoz} for details).

We            recall            very            briefly            the
construction.  Following~\cite{Avila:Gouezel:Yoccoz}   we  define  the
\emph{type}   of   $T$    by   $t$   if   $\lambda_{\pi_t^{-1}(d)}   >
\lambda_{\pi_b^{-1}(d)}$   and  $b$   if   $\lambda_{\pi_t^{-1}(d)}  <
\lambda_{\pi_b^{-1}(d)}$. When $T$ is  of type $t$ (respectively, $b$)
we   will   say   that   the  label   $\pi_t^{-1}(d)$   (respectively,
$\pi_b^{-1}(d)$) is the winner and that $\pi_b^{-1}(d)$ (respectively,
$\pi_t^{-1}(d)$) is the looser.  We define a subinterval $J$ of $I$ by
$$  J=\left\{ \begin{array}{ll}  I  \backslash T(I_{\pi_b^{-1}(d)})  &
\textrm{if  $T$ is  of  type t};\\  I  \backslash I_{\pi_t^{-1}(d)}  &
\textrm{if $T$ is of type b.}
\end{array} \right.
$$ The image of $T$  by the {\cita Rauzy-Veech} induction $\mathcal R$
is defined  as the  first return  map of $T$  to the  subinterval $J$.
This  is again  an interval  exchange transformation,  defined  on $d$
letters  (see  e.g.~\cite{Rauzy}).  The  data  of  $\mathcal R(T)$ are very easy to express 
in term of those of $T$.

There are two cases to distinguish depending whether $T$ is of type
$t$ or $b$; the labeled permutations of $\mathcal{R}(T)$ only depends on
$\pi$ and  on the type of  $T$. If $\varepsilon \in \{t,b\}$ is the type of $T$, 
this defines  two maps $\mathcal{R}_t$ and  $\mathcal{R}_b$ by $\mathcal{R}(T)=(\mathcal{R}_\varepsilon(\pi),
\lambda^{\prime})$. We
will often make use of the following notation: if $\varepsilon\in \{t,b\}$ we denote by 
$1-\varepsilon$ the other element of $\{t,b\}$. 
\begin{enumerate}
\item $T$ has type  $t$. Let $k\in \{1,\dots ,d-1\}$ such that $\pi_b^{-1}(k)=\pi_t^{-1}(d)$.  Then $\mathcal  R_t(\pi_t,\pi_b)=(\pi_t',\pi_b')$ where
$\pi_t=\pi_t'$ and
$$ \pi_b'^{-1}(j) = \left\{
\begin{array}{ll}
\pi_b^{-1}(j) & \textrm{if $j\leq  k$}\\ \pi_b^{-1}(d) & \textrm{if $j
= k+1$}\\ \pi_b^{-1}(j-1) & \textrm{otherwise.}
\end{array} \right.
$$

\item $T$ has type  $b$. Let $k\in \{1,\dots ,d-1\}$ such that $\pi_t^{-1}(k)=\pi_b^{-1}(d)$. Then $\mathcal  R_b(\pi_t,\pi_b)=(\pi_t',\pi_b')$ where
$\pi_b=\pi_b'$ and
$$ \pi_t'^{-1}(j) = \left\{
\begin{array}{ll}
\pi_t^{-1}(j) & \textrm{if $j\leq  k$}\\ \pi_t^{-1}(d) & \textrm{if $j
= k+1$}\\ \pi_t^{-1}(j-1) & \textrm{otherwise.}
\end{array} \right.
$$

\item Let  us denote  by $E_{\alpha\beta}$ the  $d\times d$  matrix of
which the  $\alpha,\beta$-th element  is equal to  $1$, all  others to
$0$.   If  $T$   is   of   type  $t$   then   let  $(\alpha,\beta)   =
(\pi^{-1}_t(d),\pi^{-1}_b(d))$   otherwise   let   $(\alpha,\beta)   =
(\pi^{-1}_b(d),\pi^{-1}_t(d))$.   Then  if  $V_{\alpha\beta}$  is  the
transvection          matrix          $I+E_{\alpha\beta}$         then
$V_{\alpha\beta}\lambda'=\lambda$.

\end{enumerate}

\begin{Remark}
In  the {\cita Veech}'s  original construction,  the matrices  used to
obtain $\lambda'$ in  terms of $\lambda$ were more  complicated: of the
form $P+E_{\alpha\beta}$  where $P$  is a permutation  matrix. Indeed,
after the {\cita Rauzy} induction ``bottom'', we usually have $\pi_t'\neq Id$, and we
must ``renumber'' it.
\end{Remark}

This construction  is due to  {\cita Rauzy}. This induction  is called
the {\cita  Rauzy-Veech} induction  since {\cita Veech}  observed that
one can  actually define the induction  on the suspension  data in the
following  way. If $\tau$  is a  suspension data  over $(\pi,\lambda)$
then we define $\mathcal R(\pi,\lambda,\tau)$ by
$$        \mathcal        R(\pi,\lambda,\tau)       =        (\mathcal
R_{\varepsilon}(\pi),V^{-1}\lambda,V^{-1}\tau),
$$ where $\varepsilon$ is the type of $T=(\pi,\lambda)$ and $V$ is the
corresponding      transition     matrix.      In      other     terms
$V_{\alpha\beta}\zeta'=\zeta$ where $\zeta=(\lambda,\tau)$.

\begin{Remark}
\label{rk:isometric}
By  construction  the  two  translation  surfaces  $X(\pi,\zeta)$  and
$X(\pi',\zeta')$ are naturally isometric (as translation surfaces).
\end{Remark}

Now  if we  iterate the  {\cita Rauzy}  induction, we  get  a sequence
$(\alpha_k,\beta_k)$          of          winners/loosers.          If
$\mathcal{R}^{(n)}(\pi,\lambda)=(\pi^{(n)},\lambda^{(n)})$     then    the
transition  matrix  that  rely  $\lambda^{(n)}$ to  $\lambda$  is  the
product of the transition matrices:
\begin{eqnarray}
\label{eq:path}
\left(\prod_{k=1}^{n} V_{\alpha_k\beta_k}\right)\lambda^{(n)}=\lambda.
\end{eqnarray}

For  a labeled  permutation $\pi$,  we call  the  \emph{labeled {\cita
Rauzy}  diagram},  denoted  by  $\mathcal  D(\pi)$,  the  graph  whose
vertices are all  labeled permutations that we can  obtained from $\pi$
by the  combinatorial {\cita Rauzy} moves.  From  each vertices, there
are two edges labeled $t$ and  $b$ (the type) corresponding to the two
combinatorial   {\cita  Rauzy}   moves.   We  will   denote  by   $\pi
\xrightarrow{\alpha,\beta}  \pi'$   for  the  edge   corresponding  to
$\mathcal R_\varepsilon(\pi) = \pi'$ where $\varepsilon\in\{t,b\}$ and
$\alpha/\beta$  is the winner/looser.  To each  path $\gamma$  is this
diagram, there is thus a sequence  of winners/loosers. We will denote by
$V(\gamma)$    the   product   of    the   transition    matrices   in
Equation~(\ref{eq:path}). The next lemma is clear from the definition.

\begin{Lemma}
\label{lm:winner:looser}
Let $\gamma_n=\pi_1\ldots\pi_n$ be a path in the labeled {\cita Rauzy}
diagram,  and  let  $V_n$  be   the  matrix  associated  to  the  path
$\gamma_n$. Let $\alpha,\beta$ be  the winner/looser associated to the
edge  $\pi_{n-1}  \xrightarrow{\alpha,\beta}  \pi_n$.  Then  $V_n$  is
obtained from  $V_{n-1}$ by  adding the column $\alpha$ to the  column $\beta$.
\end{Lemma}

\begin{Definition}
A  closed path  in the  labeled {\cita  Rauzy} diagram  is said  to be
primitive if  the associated matrix  $V$ is primitive,  \emph{i.e.} if
there exists a power of $V$ such that all the entries are positive. We
will  also say  that a  path contains  the letter  $\alpha$  as winner
(respectively,    looser)   if   it    contains   the    edge   $\cdot
\xrightarrow{\alpha,\beta}      \cdot$      (respectively,      $\cdot
\xrightarrow{\beta,\alpha} \cdot$), for some $\beta$.
\end{Definition}

We  have the  following  proposition (see~\cite{Marmi:Moussa:Yoccoz},
proposition in section $1.2.3$),

\begin{Proposition}[{\cita Marmi, Moussa \& Yoccoz}]
\label{prop:closed}
A closed path $\gamma$ in a labeled {\cita Rauzy} diagram is primitive
if and  only if $\gamma$ contains  all the letters as  winner at least
once.
\end{Proposition}

\subsection{Reduced {\cita Rauzy} diagrams}
We have  previously defined {\cita Rauzy} induction  and {\cita Rauzy}
diagrams for labeled interval exchange transformations.  One can also
define the  same for reduced interval  exchange transformations, as
it was first, for which  the corresponding labeled permutation is just
a  permutation of  $\{1,\ldots,d\}$ (see~\cite{Veech1982}).  These are
obtained after  identifying $(\pi_t,\pi_b)$ with  $(\pi_t',\pi_b')$ if
$\pi_b\circ\pi_t^{-1}=\pi_b'\circ\pi_t'^{-1}$.   In the  next  we will
use the  notation $\mathcal D_r(\pi)$  to denote the  reduced {\cita
Rauzy} diagram associated to the permutation $\pi$.

Note that the labeled {\cita Rauzy} diagram is naturally a covering of
the  reduced  {\cita  Rauzy}  diagram,  and  they  are  usually  not
isomorphic.

Given a closed  path $\gamma$ in the reduced  {\cita Rauzy} diagram,
as  previously, one  can associate  a matrix  $V$ as  follow:  we take
$(\pi_t,\pi_b)$ the labeled  permutation corresponding to the endpoint
of $\gamma$ so that $\pi_t=Id$. Then we consider $\hat{\gamma}$ a lift
of  $\gamma$   in  the  labeled  {\cita  Rauzy}   diagram.   The  path
$\hat{\gamma}$ is not  necessarily closed and it ends at a permutation
$(\pi_t',\pi_b')$. We can associate to it  a matrix $\widehat{V}$
as before. Let $P$ be  the permutation matrix defined by permuting the
columns  of   the  $d\times  d$  identity  matrix   according  to  the
permutation $\pi_t'$, \emph{i.e.} the $P=[p_{ij}]$, with $p_{ij}=1$ if
$j=\pi'_t(i)$ and  $0$ otherwise. The transition  matrix associated to
the path is then:
\begin{eqnarray}
\label{eq:matrix:unlabel}
V = \widehat{V} \cdot P.
\end{eqnarray}
As before,  a closed  path in the  reduced {\cita Rauzy}  diagram is called 
\emph{primitive} if $V$ is primitive.
A standard reference for the next two sections is~\cite{Veech1982}.

\subsection{Construction of pseudo-{\cita Anosov} homeomorphisms}

There  is a  natural  $\textrm{SL}_{2}(\R)$-action on  the strata.  In
particular, the one parameter subgroup $g_t= \left(
\begin{smallmatrix}  e^t &  0    \\  0  &    e^{-t}
\end{smallmatrix}\right)$ is called the {\cita Teichm\"uller} geodesic
flow.   It can  be  shown that  conjugacy  classes of  pseudo-{\cita
Anosov}  homeomorphisms are one-to-one with closed  geodesics of
the {\cita Teichm\"uller} geodesic flow  on strata. There is a very nice
construction of pseudo-{\cita  Anosov} homeomorphisms using the {\cita
Rauzy-Veech} induction, and we recall now this construction. \medskip

Let $\pi$ be  an irreducible permutation and let  $\gamma$ be a closed
loop in the  reduced {\cita Rauzy} diagram associated  to $\pi$. One
can associate  to $\gamma$ a  matrix $V(\gamma)$ (see  section above).
Let us assume that $V$ is {\it  primitive} and let $\theta > 1$ be its
{\cita Perron-Frobenius} eigenvalue.  We choose a positive eigenvector
$\lambda$   for   $\theta$.   It    can   be   shown   that   $V$   is
symplectic~\cite{Veech1982}, thus let  us choose an eigenvector $\tau$
for the  eigenvalue $\theta^{-1}$ with  $\tau_{\pi^{-1}_{0}(d)}>0$. It
turns   out    that   $\tau$   defines   a    suspension   data   over
$T=(\pi,\lambda)$. Indeed, the set of suspension data is an open cone,
that  is  preserved by  $V^{-1}$.  Since  the  matrix $V^{-1}$  has  a
dominant eigenvalue $\theta$ (for  the eigenvector $\tau$), the vector
$\tau$ must belong to this cone. If $\zeta=(\lambda,\tau)$, one has
\begin{multline*}
\mathcal{R}(\pi,\zeta) = (\pi,V^{-1}\zeta) = (\pi,V^{-1}\lambda,V^{-1}
\tau)   =   (\pi,\theta^{-1}\lambda,\theta   \tau)   =  \\   =   g_{t}
(\pi,\lambda,\tau), \qquad \textrm{ where } \qquad t=\log(\theta) > 0.
\end{multline*}
Hence the  two surfaces $X(\pi,\zeta)$  and $g_{t}X(\pi,\zeta)$ differ
by    some    element    of    the   mapping    class    group    (see
Remark~\ref{rk:isometric}).    In   other   words   there   exists   a
pseudo-{\cita  Anosov}  homeomorphism   $\pA$ affine with  respect  to  the
translation surface  $X(\pi,\zeta)$ and such that $D \pA  = g_{t}$. The action 
of $\pA$ on the relative homology of $(X,\omega)$ is $V(\gamma)$ thus 
the dilatation of  $\pA$ is $\theta$. Note that  by construction $\pA$
fixes the zero on the left of the interval $I$ and also a horizontal separatrix
adjacent to this zero (namely, the oriented half line corresponding to the interval $I$).

It turns out that this construction is very general as we  will see in the
coming section.

\subsection{Discrete representation of the geodesic flow}
\label{subsec:Veech}

Let  us  fix an  irreducible  permutation  $\pi$  defined over $d$  letters.  If
$\mathcal  C$   is  the  connected  component  of   some  stratum,  let
$\widehat{\mathcal  C}$  be  the  ramified  cover  over  $\mathcal  C$
obtained  by  considering the  set  of  triplets $(X,\omega,l)$  where
$(X,\omega) \in \mathcal C$ and $l$ is a horizontal separatrix adjacent to a zero of $\omega$. \medskip

Clearly the set  of $(\lambda,\tau)$ such that $\tau$  is a suspension
data over $\pi$ is a connected space and the map $(\lambda,\tau)=\zeta
\mapsto  X(\pi,\zeta)$ is  continuous. Thus  all surfaces  obtained by
this  construction belong  to  the same  connected  component of  some
strata, say $\mathcal C(\pi) \subset \HA(\sigma)$. Moreover
$\sigma$ can be computed easily in terms of $\pi$.  We define
$$   \mathcal   T(\widehat{\mathcal   C})  =   \left\{   (\pi,\zeta);\
\widehat{\mathcal C}(\pi) = \widehat{\mathcal C}, \textrm{ and } \zeta
\textrm{ is a suspension datum for } \pi \right\}.
$$  The  {\cita Rauzy-Veech}  induction  is  (almost everywhere)  well
defined  and one-to-one on  $\mathcal A(\widehat{\mathcal  C})$. Hence
let $\mathcal  H(\widehat{\mathcal C})$  be the quotient  of $\mathcal
T(\widehat{\mathcal C})$ by the induction.

The {\cita  Veech} zippered rectangle's  construction provides (almost
everywhere)  a one-to-one  map $Z  : \mathcal  H(\widehat{\mathcal C})
\rightarrow \widehat{\mathcal C} $ (see~\cite{Boissy:rc} for details).
$$
\begin{array}{ccccc}
\HA(\widehat{\mathcal    C})    &    \overset{Z}{\longrightarrow}    &
\widehat{\mathcal  C} \\  &&  \downarrow\\ &&\mathcal  C  & \subset  &
\HA(\sigma)
\end{array} 
$$  One  can  define   the  {\cita  Teichm\"uller}  geodesic  flow  on
$\HA(\widehat{\mathcal C})$ by $g_t (\pi,\zeta)= (\pi,g_t \zeta)$. The
{\cita Teichm\"uller} flow on $\mathcal  C$ lifts to a flow $g_{t}$ on
$\widehat{\mathcal C}$. It is easy  to check that $g_t$ is equivariant
with $Z$ i.e. $g_{t} Z = Z g_{t}$. \medskip

By  construction,   periodic  orbits   of  $g_{t}$  on   $\mathcal  C$
corresponding  to pseudo-{\cita  Anosov} homeomorphisms that  fix the  separatrix $I$
lift to  periodic orbits on  $\widehat{\mathcal C}$ for  $g_{t}$. Thus
they    produce   periodic    orbits   for    $g_{t}$    on   suspensions
$\HA(\widehat{\mathcal C})$. \medskip

\noindent In  fact   {\cita  Veech}  proved   more: all pseudo-{\cita  Anosov} homeomorphisms 
fixing a separatrix arise in this way. The  subset   of  $\mathcal
T(\widehat{\mathcal C})$ defined by
$$ \left\{(\pi,\zeta)\in \mathcal{T}(\widehat{\mathcal  C});\ 1 \leq |
Re(\zeta)|                                                         \leq
1+\min\bigl(Re(\zeta_{\pi_{0}^{-1}(d)}),Re(\zeta_{\pi_1^{-1}(d)})
\bigr) \right\}
$$ is  a fundamental domain of $\mathcal  T(\widehat{\mathcal C})$ for
the   quotient  map  $\mathcal{T}(\widehat{\mathcal   C})  \rightarrow
\mathcal{H}(\widehat{\mathcal C})$, and  the {\cita Poincar\'e} map of
the {\cita Teichm\"uller} flow on the section
$$   \mathcal{S}=\{(\pi,\zeta);\   \pi   \textrm{   irreducible},\   |
Re(\zeta)|=1\}/ \sim
$$ is  precisely the {\it renormalized}  {\cita Rauzy-Veech} induction
on suspensions:
$$           \widehat{\mathcal{R}}(\pi,\lambda):=           (\mathcal
R_{\varepsilon}(\pi),V^{-1}\lambda/|V^{-1}\lambda|).
$$

We can summarize the above discussion by the following theorem.

\begin{NoNumberTheorem}[{\cita Veech}]
Let $\gamma$ be  a closed loop, based at~$\pi$,  in a reduced {\cita
Rauzy}  diagram $\mathcal  D_{r}(\pi)$  and let  $V=V(\gamma)$ be  the
product of the associated transition  matrices. Let us assume that $V$
is primitive.  Let $\lambda$ be  a positive eigenvector for the {\cita
Perron-Frobenius}  eigenvalue  $\dil$ of  $V$  and  let  $\tau$ be  an
eigenvector   for    the   eigenvalue   $\dil^{-1}$    of   $V$   with
$\tau_{\pi_0^{-1}(1)} > 0$. Then
\begin{enumerate}

\item    $\zeta=(\lambda,\tau)$   is    a    suspension   datum    for
$T=(\pi,\lambda)$;

\item  The matrix $A=\left(  \begin{smallmatrix}\dil^{-1} &  0 \\  0 &
\dil  \end{smallmatrix} \right)$ is  the derivative  map of  an affine
pseudo-{\cita  Anosov}  homeomorphism  $\pA$  on  $X(\pi,\zeta)$;  The
action on relative homology of $\pA$ is the matrix $V$.

\item The dilatation of $\phi$ is $\dil$;

\item All  pseudo-{\cita Anosov} homeomorphisms that  fix a separatrix
are obtained in this way.
\end{enumerate}

\end{NoNumberTheorem}

We will use this theorem in order to prove our main result.

\subsection{Examples of labeled and reduced {\cita Rauzy} diagrams}
\label{subsec:hyper}

\begin{Convention}
Let $\pi$ be  a permutation. We will denote  by $\mathcal{D}(\pi)$ the
labeled {\cita Rauzy} diagram of $\pi$ and by $\mathcal{D}_r(\pi)$ the
reduced one.
\end{Convention}

\subsubsection{Hyperelliptic connected components}
\label{sec:description:unmarked}
Let  $n\geq  2$.  A   representative  permutation  for  the  connected
component    $\mathcal    H^{hyp}(2g-2)$   (respectively,    $\mathcal
H^{hyp}(g-1,g-1)$) is
\begin{eqnarray}
\label{eq:permutation:hyp}
\tau_n = \left( \begin{smallmatrix} 0 & 2 & 3 & \dots & n-1 & n \\ n &
n-1 & \dots & 3 & 2 & 0 \end{smallmatrix} \right),
\end{eqnarray}
where $n=2g$ (respectively, $n=2g+1$). We use this permutation in order to simplify 
the notations latter.  It turns out that labeled and
reduced   {\cita  Rauzy}  diagrams   are  isomorphic.   The  precise
description    of     the    diagrams    was     given    by    {\cita
Rauzy}~\cite{Rauzy}. Let us recall the result here. \medskip

If $\pi=(\pi_t,\pi_b)$  is a labeled  permutation, for $\varepsilon\in
\{t,b\}$ we  define $G^\varepsilon(\pi)$ to  be the subdiagram  of the
{\cita Rauzy} diagram of $\pi$  whose vertices are obtained from $\pi$
by   a   simple   path,   and    whose   first   step   is   the   map
$\mathcal{R}_\varepsilon$.  Recall that for  $\varepsilon\in \{t,b\}$,
we denote by $1-\varepsilon$ the other element of $\{t,b\}$.
\begin{NoNumberProposition}[{\cita Rauzy}]
Let    $\tau_n$    be    the    labeled   permutation    defined    in
Equation~(\ref{eq:permutation:hyp}). Then
\begin{enumerate}
\item The  vertices of  the {\cita Rauzy}  diagram of $\tau_n$  is the
disjoint union  of the  vertices of $G^0(\tau_n)$,  $G^1(\tau_n)$, and
$\{\tau_n\}$.
\item  Let  $\tau_{n,k}=\mathcal{R}_\varepsilon^k(\tau_n)$,  for  some
$\varepsilon\in    \{t,b\}$   and   $k\in    \{1,\ldots,n-1\}$.   Then
$G^{1-\varepsilon}(\tau_{n,k})$    is    naturally    isomorphic    to
$G^{1-\varepsilon}(\tau_{n-k})$.
\item The cardinality of the {\cita Rauzy} diagram is $2^{n-1}-1$.
\end{enumerate}
\end{NoNumberProposition}

\begin{figure}
\begin{center}
\psfrag{1}[]{{\tiny        $\left(\begin{smallmatrix}        0&2&3&4\\
4&3&2&0\end{smallmatrix}\right)$}}                   \psfrag{2}[]{\tiny
$\left(\begin{smallmatrix} 0&2&3&4\\ 4&0&3&2\end{smallmatrix}\right)$}
\psfrag{3}[]{\tiny         $\left(\begin{smallmatrix}        0&2&4&3\\
4&2&0&3\end{smallmatrix}\right)$}                    \psfrag{4}[]{\tiny
$\left(\begin{smallmatrix} 0&4&2&3\\ 4&3&0&2\end{smallmatrix}\right)$}
\psfrag{5}{\tiny          $\left(\begin{smallmatrix}         0&3&4&2\\
4&3&0&2\end{smallmatrix}\right)$}                    \psfrag{6}[]{\tiny
$\left(\begin{smallmatrix} 0&2&4&3\\ 4&0&3&2\end{smallmatrix}\right)$}
\psfrag{7}[]{\tiny         $\left(\begin{smallmatrix}        0&4&2&3\\
4&3&2&0\end{smallmatrix}\right)$} \includegraphics[scale=0.6]{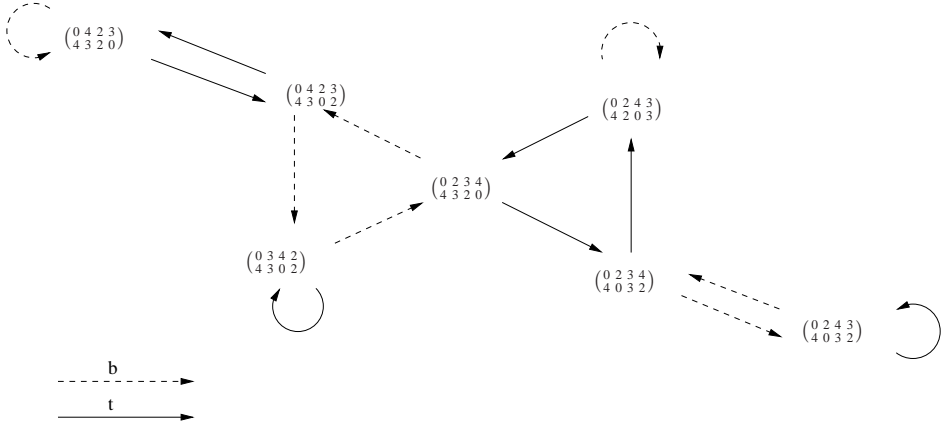}
\caption{The diagram $\mathcal D^{hyp}  = \mathcal  D(\tau_{n})  \approx  \mathcal  D_r(\tau_{n})$ for $n=4$.}
\end{center}
\end{figure}

Using  this  result, one  can  show  that  the labeled  and  reduced
diagrams  are  isomorphic.  We  will  denote  by  $\mathcal D^{hyp}  =
\mathcal  D(\tau_{n})  \approx  \mathcal  D_r(\tau_{n})$  this  {\cita
Rauzy}  diagram.  A  consequence of  the previous  proposition  is the
following:
\begin{Corollary} \label{Cor:Dhyp}
Let $\gamma$ be a closed (oriented) path in $\mathcal D^{hyp}$. Assume
that it  contains the step $\pi  \xrightarrow{\alpha,\beta} \pi'$ with
$\pi\neq   \pi'$,    then   it    must   contain   the    step   $\pi'
\xrightarrow{\alpha,\beta'} \pi''$.
\end{Corollary}

\begin{proof}
A consequence  of the {\cita Rauzy}  description of $\mathcal  D^{hyp}$ is
that  $\mathcal  D^{hyp} \backslash  \{\pi'\}$  is  not connected  any
more.     Assume    that,     in    the     step     following    $\pi
\xrightarrow{\alpha,\beta} \pi'$, the  symbol $\alpha$ is looser, then
the new permutation and $\pi$ belong to two different connected components
of  $\mathcal D^{hyp}  \backslash \pi'$.  Looking at  the  sequence of
permutations  that  appear,  one  must  come  back  to  the  connected
component     of    $\pi$,     and    hence,     the     step    $\pi'
\xrightarrow{\alpha,\beta'} \pi''$ eventually appears.
\end{proof}

\subsubsection{Hyperelliptic connected components with a marked point}
\label{sec:description:marked}

Here  we present  a second  family  of {\cita  Rauzy} diagrams that will be useful for our proof.  There
description  is  a little  bit  more  technical,  but will  be  needed
later. We start with a very informal description of the labeled {\cita
Rauzy}  diagram. The  precise description  can be  skipped in  a first
reading. \medskip

The  labeled {\cita  Rauzy} diagram  is  a covering  of the  reduced
one. The  cardinalities of these  labeled {\cita Rauzy}  diagrams have
been  calculated  (see  {\cita Delecroix}~\cite{Delecroix2010}).   The
degree of the covering is also known~\cite{Boissy:degree}. \medskip

\begin{figure}
\begin{center}
\psfrag{1a}[]{{\tiny       $\left(\begin{smallmatrix}      0&2&3&1&4\\
4&3&2&1&0\end{smallmatrix}\right)$}}                \psfrag{2a}[]{\tiny
$\left(\begin{smallmatrix}                                  0&2&3&1&4\\
4&0&3&2&1\end{smallmatrix}\right)$}                 \psfrag{8a}[]{\tiny
$\left(\begin{smallmatrix}                                  0&4&2&3&1\\
4&3&2&1&0\end{smallmatrix}\right)$}                 \psfrag{5a}[]{\tiny
$\left(\begin{smallmatrix}                                  0&2&4&3&1\\
4&1&0&3&2\end{smallmatrix}\right)$}                   \psfrag{5c}{\tiny
$\left(\begin{smallmatrix}                                  0&1&4&2&3\\
4&3&0&2&1\end{smallmatrix}\right)$}                 \psfrag{7a}[]{\tiny
$\left(\begin{smallmatrix}                                  0&2&3&4&1\\
4&2&1&0&3\end{smallmatrix}\right)$}                 \psfrag{7b}[]{\tiny
$\left(\begin{smallmatrix}                                  0&3&1&4&2\\
4&3&2&0&1\end{smallmatrix}\right)$}  \psfrag{3a}[]{*} \psfrag{4a}[]{*}
\psfrag{6a}[]{*}  \psfrag{9a}[]{*} \psfrag{10a}[]{*} \psfrag{12a}[]{*}
\includegraphics[scale=0.6]{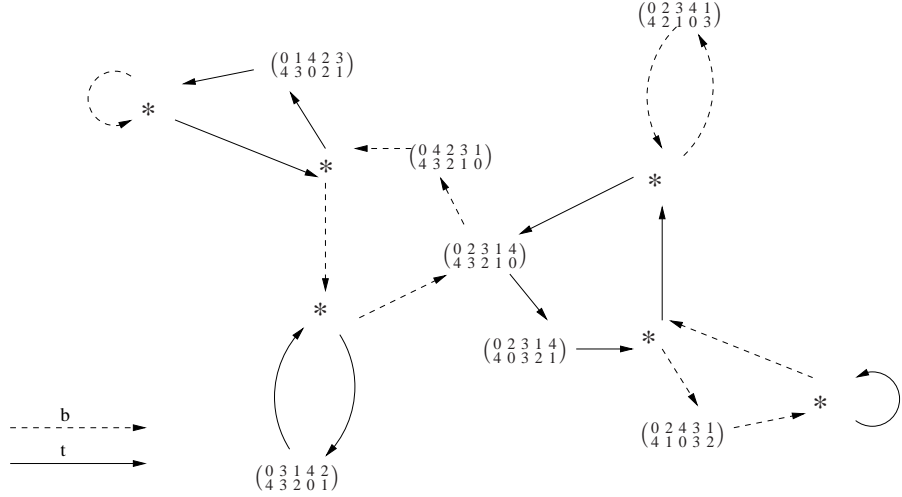}
\caption{The augmented  hyperelliptic diagram, and  added permutations
for $n=4$.}
\label{augmented:hyp}
\end{center}
\end{figure}

\begin{figure}[htbp]
\begin{center}
\psfrag{1a}[]{{\tiny       $\left(\begin{smallmatrix}      0&2&3&1&4\\
4&3&2&1&0\end{smallmatrix}\right)$}}                \psfrag{2a}[]{\tiny
$\left(\begin{smallmatrix}                                  0&2&3&1&4\\
4&0&3&2&1\end{smallmatrix}\right)$}                 \psfrag{8a}[]{\tiny
$\left(\begin{smallmatrix}                                  0&4&2&3&1\\
4&3&2&1&0\end{smallmatrix}\right)$}  \psfrag{5a}[]{C}  \psfrag{C}{C  =
$\left(\begin{smallmatrix}                                  0&2&4&3&1\\
4&1&0&3&2\end{smallmatrix}\right)$}   \psfrag{5c}{A}   \psfrag{A}{A  =
$\left(\begin{smallmatrix}                                  0&1&4&2&3\\
4&3&0&2&1\end{smallmatrix}\right)$}

\psfrag{7a}[]{B} \psfrag{B}{B = $\left(\begin{smallmatrix} 0&2&3&4&1\\
4&2&1&0&3\end{smallmatrix}\right)$}

\psfrag{7b}[]{\tiny       $\left(\begin{smallmatrix}       0&3&1&4&2\\
4&3&2&0&1\end{smallmatrix}\right)$}  \psfrag{3a}[]{*} \psfrag{4a}[]{*}
\psfrag{6a}[]{*}  \psfrag{9a}[]{*} \psfrag{10a}[]{*} \psfrag{12a}[]{*}
\psfrag{1b}[]{{\tiny $\left(\begin{smallmatrix} 0&3&1&2&4\\ 4&1&3&2&0
\end{smallmatrix}\right)$}} \psfrag{2b}[]{\tiny
$\left(\begin{smallmatrix}                                  0&3&1&2&4\\
4&0&1&3&2\end{smallmatrix}\right)$}                 \psfrag{8b}[]{\tiny
$\left(\begin{smallmatrix}                                  0&4&3&1&2\\
4&1&3&2&0\end{smallmatrix}\right)$}  \psfrag{5b}[]{D}  \psfrag{D}{D  =
$\left(\begin{smallmatrix}                                  0&3&4&1&2\\
4&2&0&1&3\end{smallmatrix}\right)$}                 \psfrag{7c}[]{\tiny
$\left(\begin{smallmatrix}                                  0&1&2&4&3\\
4&1&3&0&2\end{smallmatrix}\right)$}  \psfrag{3b}[]{*} \psfrag{4b}[]{*}
\psfrag{6b}[]{*}  \psfrag{9b}[]{*} \psfrag{10b}[]{*} \psfrag{12b}[]{*}
\psfrag{1c}[]{{\tiny $\left(\begin{smallmatrix} 0&1&2&3&4\\ 4&2&1&3&0
\end{smallmatrix}\right)$}} \psfrag{2c}[]{\tiny
$\left(\begin{smallmatrix}                                  0&1&2&3&4\\
4&0&2&1&3\end{smallmatrix}\right)$}                 \psfrag{8c}[]{\tiny
$\left(\begin{smallmatrix}                                  0&4&1&2&3\\
4&2&1&3&0\end{smallmatrix}\right)$}  \psfrag{3c}[]{*} \psfrag{4c}[]{*}
\psfrag{6c}[]{*} \psfrag{9c}[]{*} \psfrag{10c}[]{*} \psfrag{12c}[]{*}

\includegraphics[scale=0.5]{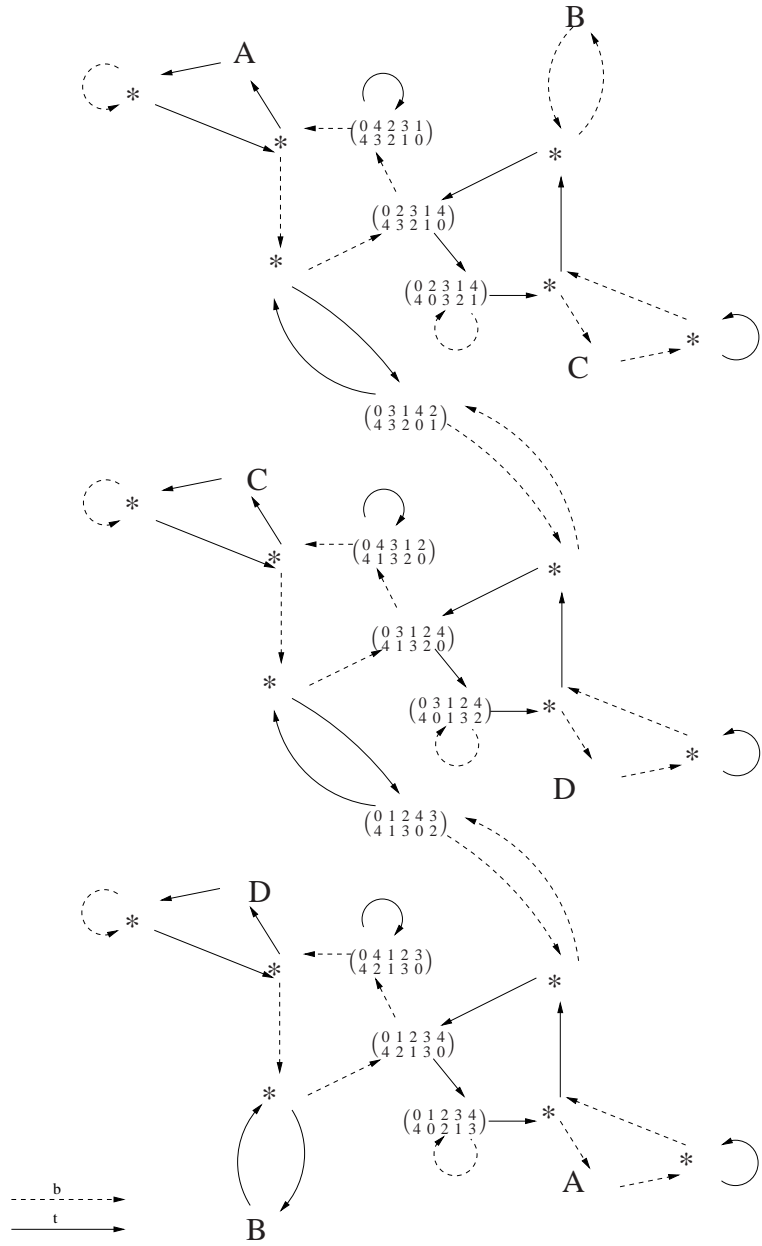}
\caption{The  complete diagram  $\mathcal{D}(\pi_n)$  for $n=4$.   The
vertices labeled by letters $A,B,C$ and $D$ should be identified.}
\label{fig:diagram:labeled}
\end{center}
\end{figure}

A fundamental domain of this covering can be roughly seen as a copy of
the  hyperelliptic {\cita  Rauzy} diagram  described  previously, with
some   added  permutations   (see   Figure~\ref{augmented:hyp}).  This
fundamental domain (a \emph{leaf}) is  composed at first glance by two
principal loops that intersect in a \emph{central permutation}, and on
the over  vertices of these loops  starts a secondary  loop. The whole
diagram is  obtained by  taking several copies  of this  ``leaf''. The
different   leaves  are  joined   together  by   the  \emph{transition
permutations}: each secondary loops  contains a unique such transition
permutation, and a ``k-th secondary loop''  of a leaf is attached to a
``k-th     secondary      loop''     of     another      leaf     (see
Figure~\ref{fig:diagram:labeled},     and      compare     it     with
Figure~\ref{augmented:hyp}).
\medskip

We now give the  precise description of this diagram. A representative permutation for
the  connected  component  $\mathcal  H^{hyp}(0,2g-2)$  (respectively,
$\mathcal H^{hyp}(0,g-1,g-1)$) is
\begin{eqnarray}
\label{eq:permutation:marqued:hyp}
\pi_n = \left( \begin{smallmatrix} 0 & 2 &  3 & \dots & n-1 & 1 & n \\
n & n-1 & \dots & 3 & 2 & 1 & 0 \end{smallmatrix} \right),
\end{eqnarray}
where $n=2g$ (respectively, $n=2g+1$).  Contrary to the previous case,
the diagrams labeled  and reduced are not isomorphic  anymore. As we
will see the cardinality of the  reduced diagram is $2^{n-1} -1 + n$
and the cardinality of the labeled diagram is $(2^{n-1} -1+ n)(n-1)$.

Our  next goal  is to  describe the  reduced {\cita  Rauzy} diagrams
$\mathcal{D}(\pi_n)$ and $\mathcal D_r(\pi_n)$ associated to the above
permutation $\pi_n$. The key point is to observe that if we forbid the
letter  $1$ to be  winner or  looser in the  construction of  the diagram
$\mathcal D(\pi_n)$  starting from $\pi$,  then one gets more  or less
the  diagram   $\mathcal  D^{hyp}$  (compare  with   {\cita  Avila  \&
Viana}~\cite{Avila:Viana}). So  we have to  add the letter $1$  to the
permutations in $\mathcal D^{hyp}$. One can do that as follows. \medskip

\begin{enumerate}
\item To each  permutation $\pi \in \mathcal D^{hyp}$  construct a new
permutation $\tilde  \pi \in \mathcal  D(\pi_n)$ by adding  the letter
$1$ to the  left of the letter $n$  on the top and to the  left of the
letter  $0$  on  the bottom.  We  will  refer  to the  {\it  augmented
hyperelliptic diagram}.
\item
\label {eq:add:vertex:2}
Replace each  edge $\pi \xrightarrow{k,n} \pi'$  in $\mathcal D^{hyp}$
by  two edges  $\tilde \pi  \xrightarrow{k,n}  \pi'' \xrightarrow{k,1}
\tilde \pi'$ in $\mathcal D(\pi_n)$, where the permutation $\pi''$ has
the letter $1$ for the last letter on the top.
\item
\label {eq:add:vertex:3}
Replace each  edge $\pi \xrightarrow{k,0} \pi'$  in $\mathcal D^{hyp}$
by  two edges  $\tilde \pi  \xrightarrow{k,0}  \pi'' \xrightarrow{k,1}
\tilde \pi'$ in $\mathcal D(\pi_n)$, where the permutation $\pi''$ has
the letter $1$ for the last letter on the bottom.

\item We will  denote by $\mathcal A_n \subset  \mathcal D(\pi_n)$ the
added    permutations    of    the    operations~(\ref{eq:add:vertex:2})
and~(\ref{eq:add:vertex:3}). All  the edges  are built except  the cases
when $\pi'' \in \mathcal A_n $ and where $1$ is winner.
\end{enumerate}

Observe              that              operation~(\ref{eq:add:vertex:2})
(respectively,~(\ref{eq:add:vertex:3}))  arises exactly when  the letter
$n$   (respectively,  $0$)   is  looser   in  the   diagram  $\mathcal
D^{hyp}$. Hence the description given  in the previous section of this
hyperelliptic diagram leads to the following.

\begin{Lemma}
\label{lm:description}
In  $\mathcal{D}^{hyp}$, the  edges where  $n$ (respectively,  $0$) is
looser are
$$    \mathcal{R}^k_t    (\pi)    \xrightarrow{*,n}    \mathcal    R_b
\mathcal{R}^k_t (\pi)  \qquad \textrm{(respectively, } \mathcal{R}^k_b
(\pi) \xrightarrow{*,0} \mathcal R_t \mathcal{R}^k_b (\pi)),
$$ for  $k\in \{0,\ldots,n-2\}$. So,  there are there are  $n-1$ edges
where $n$ is looser and $n-1$ edges where $0$ is looser.

Therefore,  the  added  permutations  in $\mathcal  D(\pi_n)$  by  the
operations (2) and (3) are:
\begin{eqnarray}
\label{eq:perm1}
\mathcal R_b \mathcal{R}^k_t (\pi_n) &=&\left( \begin{smallmatrix} 0 &
2 & \dots&k-1& k&n&k+1&\ldots & n-1 & 1 \\ n & k-1 &\ldots & 2 & 1 & 0
& n-1 & \dots & k+1&k \end{smallmatrix} \right), \\
\label{eq:perm2}
\mathcal  R_t  \mathcal{R}^k_b  (\pi_n) &=&\left(  \begin{smallmatrix}
0&n-k+2&\ldots&n-1&1&n&2&\ldots&n-k+1 \\ n & n-1 &\ldots&n-k+2 &n-k+1&
0& n-k &\ldots& 1 \end{smallmatrix} \right),
\end{eqnarray}
for $k\in \{2,\ldots,n-1\}$, and
\begin{eqnarray}
\label{eq:perm3}
\mathcal{R}_t(\pi_n)&=&\left( \begin{smallmatrix} 0 &  2 & 3 & \dots &
n-1 &  1 &  n \\ n  &0 &  n-1 & \dots  & 3 &  2 &  1 \end{smallmatrix}
\right) \\
\label{eq:perm4}
\mathcal{R}_b(\pi_n)&=&\left( \begin{smallmatrix} 0&n &  2 & 3 & \dots
& n-1 & 1 \\ n & n-1 & \dots & 3 & 2 & 1 & 0 \end{smallmatrix} \right)
\end{eqnarray}
In addition  for all $k=2,\dots,n-1$  the two new permutations  in the
labeled   diagram  $\mathcal  D(\pi_n)$   defined  by~(\ref{eq:perm1})
and~(\ref{eq:perm2})  correspond  to   the  same  permutation  in  the
reduced diagram  $\mathcal D_r(\pi_n)$; the  renumbering corresponds
to    $\sigma^{-k}$    with    $\sigma$   the    cyclic    permutation
$(1,2,\ldots,n-1)$.   The  two   new  permutations   corresponding  to
$\mathcal  R_t(\pi_n)$  and  $\mathcal  R_b(\pi_n)$ are  different  in
$\mathcal D_{r}(\pi_{n})$.
\end{Lemma}

\begin{proof}[Proof of Lemma~\ref{lm:description}]
The proof is obtained by  straightforward computation and is left to
the reader.
\end{proof}

Now in  order to  finish the construction  of $\mathcal  D(\pi_n)$ and
$\mathcal{D}_r(\pi_n)$,   one   has  to   consider   from  the   added
permutations $\mathcal A_n$ the operations top/bottom where the letter
$1$  is  winner. In  fact  it turns  out  that  the reduced  diagram
$\mathcal D_{r}(\pi_n)$ is already constructed. Namely,

\begin{Corollary}
\label{cor:unlabel}
The  diagram $\mathcal  D_{r}(\pi_n)$  corresponds to  adding the  new
permutations  defined in Lemma~\ref{lm:description},  namely $\mathcal
T_n$ to  the augmented hyperelliptic  diagram, up to  renumbering.  In
particular the cardinality of this diagram is $2^{n-1}-1 + n$.
\end{Corollary}

\begin{proof}[Proof of Corollary~\ref{cor:unlabel}]
Let $\pi' \in \mathcal A_n $  be an added permutation.  Assume that it
is   given   by   (\ref{eq:perm1}),   that   is   $\pi'=\mathcal   R_b
\mathcal{R}_t^k(\pi_n)$  for  some  $k$.  The  edge  $\pi'\rightarrow
\mathcal{R}_b(\pi')$ was constructed when defining $\mathcal{T}_n$. By
Lemma~\ref{lm:description},  we  also  have  that  $\pi'=\mathcal  R_t
\mathcal{R}_b^k(\pi_n)$ in the  reduced {\cita Rauzy} diagram, hence
the  edge  $\pi'\rightarrow  \mathcal{R}_t(\pi')$  was  also  already
constructed. This corresponds to $n-2$ added permutations.

The remaining  cases are when  $\pi'$ is given by  (\ref{eq:perm3}) or
(\ref{eq:perm4}). Then is clear that the array corresponding to $1$
being winner are  arrays from $\pi'$ to itself.  This corresponds to 2
permutations.

Hence, the diagram  $\mathcal{D}_r(\pi_n)$ is completely built.  Since
$\#  \mathcal D^{hyp}=2^{n-1}-1$ one  has $\#  \mathcal D_{r}(\pi_{n})
=2^{n-1} -1 + n$.

\end{proof}

Lemma~\ref{lm:description} and Corollary~\ref{cor:unlabel} imply the
following description of the diagram $\mathcal D(\pi_n)$.

\begin{Proposition}
The  augmented hyperelliptic  diagram together  with  the permutations
$\mathcal  T_n$  form  a  fundamental  domain  for  the  covering  map
$\mathcal  D(\pi_n) \longrightarrow  \mathcal D_{r}(\pi_n)$.   A fiber
for  this  covering  consists  of  labeled permutations  of  the  kind
$\{\pi,\pi\circ   \sigma,  \ldots,  \pi\circ   \sigma^{n-2}\}$,  where
$\pi\circ  \sigma^k$   means  the  labeled   permutation  $\pi$  after
renumbering with  $\sigma^k$.  In  particular the cardinality  of this
diagram is $(2^n-1 + n)(n-1)$.
\end{Proposition}

The following technical  lemma is similar  to Corollary \ref{Cor:Dhyp}. It  is a
very important lemma since it  will allow us to give information about
\emph{all} irreducible paths.
\begin{Lemma} \label{subpath}
Let $\gamma$ be a  closed oriented path in $\mathcal{D}({\pi_n})$: Let
$k\in   \{2,\ldots,n-2\}$   and   let   $k'\in  \{2,n-k\}$   and   let
$\varepsilon\in \{t,b\}$.

\begin{itemize}
\item[a)]     If    $\gamma$     contains    the     step    $\mathcal
R_\varepsilon^{k'-1}\mathcal  R_{1-\varepsilon}^{k}(\pi_n) \rightarrow
\mathcal   R_\varepsilon^{k'}\mathcal   R_{1-\varepsilon}^{k}(\pi_n)$,
then  it also contains  the step  $\mathcal R_\varepsilon^{k'}\mathcal
R_{1-\varepsilon}^{k}        (\pi_n)        \rightarrow       \mathcal
R_\varepsilon^{k'+1}\mathcal R_{1-\varepsilon}^{k}(\pi_n)$.
\item[b)]     If    $\gamma$     contains    the     step    $\mathcal
R_\varepsilon^{k-1}(\pi_n)             \rightarrow            \mathcal
R_\varepsilon^{k}(\pi_n)$,  then it also  contains the  step $\mathcal
R_\varepsilon^{k}         (\pi_n)         \rightarrow         \mathcal
R_\varepsilon^{k+1}(\pi_n)$,           or           the           step
$\mathcal{R}_{1-\varepsilon}^k(\pi_n\circ     \sigma^i)    \rightarrow
\mathcal{R}_{1-\varepsilon}^{k+1}(\pi_n  \circ  \sigma^i)$,  for  some
$i\in \{0,\ldots,n-2\}$.
\end{itemize}
This property is  also true if $\gamma$ is a  nonclosed path that starts and ends in the set $\{\pi_n, \pi_n\circ\sigma, \ldots, \pi_n\circ\sigma^{n-2}\}$.
\end{Lemma}

\begin{proof}
a)    Using   the    construction    of   $\mathcal{D}(\pi_n)$    from
$\mathcal{D}^{hyp}$, we see that,  as in Corollary \ref{Cor:Dhyp}, the
diagram            $\mathcal{D}(\pi_n)\backslash            \{\mathcal
R_\varepsilon^{k'}\mathcal    R_{1-\varepsilon}^{k}(\pi_n)    \}$   is
nonconnected.

b)  Here,  we  have  that  $\mathcal{D}(\pi_n)\backslash  \{  \mathcal
R_\varepsilon^{k}  (\pi_n), \ \mathcal{R}_{1-\varepsilon}^k(\pi_n\circ
\sigma^i)  \}$  is  nonconnected,  for  the parameter  $i$  such  that
$\mathcal{R}_{1-\varepsilon}                  \mathcal{R}_\varepsilon^k
(\pi_n)=\mathcal{R}_\varepsilon
\mathcal{R}_{1-\varepsilon}^k(\pi_n\circ \sigma^i)$.
\end{proof}

\section{Minimal dilatations in a hyperelliptic connected component}

\subsection{A key proposition}
{\cita Veech}'s  construction can  only build a  pseudo-{\cita Anosov}
homeomorphism   fixing  a   separatrix.    But  there   can  be   many
pseudo-{\cita   Anosov}   homeomorphism    that   do   not   fix   any
separatrix~\cite{Los2009,Lanneau:fixed}.  In our case, we are saved by
the following key proposition.

\begin{Proposition}
\label{theo:reduc1}
Let  $g\geq  1$  and  let  $(X,\omega)$  be a  flat  surface  in  some
hyperelliptic   connected   component   $\mathcal  H^{hyp}(2g-2)$   or
$\mathcal  H^{hyp}(g-1,g-1)$.  Then   for  any  pseudo-{\cita  Anosov}
homeomorphism $\pA$ affine with respect to $(X,\omega)$, $\pA^2$ has a
fixed  point  of negative  index  i.e. $\pA^2$  fixes  a  point and  a
outgoing separatrix issued from that point.
\end{Proposition}

\begin{proof}[Proof of Proposition~\ref{theo:reduc1}]
Let  $\tau$  be  the  hyperelliptic  involution  on  $X$.  Firstly  by
Proposition~\ref{prop:commutes}  $\tau$  commutes  with  $\pA$.   Thus
$\pA$  descends to a  pseudo-{\cita Anosov}  homeomorphism $f$  on the
sphere.   By assumption, $f$  has a  single zero  (of order  $2g-3$ or
$2g-2$ depending the case). Thus $f$ fixes  the zero  and induces a  disc pseudo-{\cita
Anosov} homeomorphism, say $g$.  Now by {\cita Brouwer}'s theorem, $g$
has a  fixed point on the disc.  Either this point is  inside the disc
(thus this  is a regular point) or  on the boundary. If  the last case
occur, then $f$ fixes the separatrices issued from the singularity. It
is then not hard to see  that either $\pA^2$ fixes a regular point and
the  outgoing  separatrix,  or  a  singular  point  and  the  outgoing
separatrices.  Proposition~\ref{theo:reduc1} is proven.
\end{proof}

\begin{Remark}
Another {\cita  Brouwer}'s theorem states that a  homeomorphism on the
plane with a periodic orbit has  a fixed point. With this theorem, it
is  easy  to  show  that  there actually  exists  a  regular  fixed  point  for
$\pA^2$. However, we believe that the case where the separatrix is attached to a 
singularity is useful  to present  the main  ideas without the technical difficulties of the other case.
\end{Remark}

The two next sections analyse the two cases depending wether the fixed
point is singular or regular.

\subsection{Case when the fixed separatrix is adjacent to the singularity}
In  this section,  we study  the case  when the  pseudo-{\cita Anosov}
homeomorphism fixes  a horizontal separatrix starting  from a singular
point. We thus have to consider components $\mathcal H^{hyp}(2g-2)$ and 
$\mathcal H^{hyp}(g-1,g-1)$. We first  prove a sufficient condition for  a primitive path to
have a transition matrix with a spectral radius greater than 2.

\begin{Proposition}\label{prop:win:los}
Let  $\gamma$ be  a primitive  closed  path in  some reduced  {\cita
Rauzy} diagram.  Let us assume that  there is a lift  $\hat \gamma$ of
$\gamma$ in  the labeled {\cita  Rauzy} diagram that contains  all the
letters as  loosers, or all the  letters as winners.   Then the {\cita
Perron-Frobenius} eigenvalue of the matrix $V(\gamma)$ is bounded from
below by $2$.
\end{Proposition}

\begin{proof}[Proof of the proposition]
If  $\hat  \gamma$ contains  all  the  looser  (respectively, all  the
winner) then Lemma~\ref{lm:winner:looser} implies that the minimum  of the sums
of the  columns (respectively, the lines) of  the corresponding matrix
$\widehat{V}$, before remuneration, is  at least $2$. Since $V(\gamma)
=  \widehat{V}  \cdot  P$  where  $P$ is  a  permutation  matrix  (see
Equation~\ref{eq:matrix:unlabel}) the same property on the sums of the
columns (respectively, the lines) holds for $V(\gamma)$.

By assumption, the matrix  $V=V(\gamma)$ is primitive. Up to replacing
$V$ by  $^tV$, which  does not change  its eigenvalues, we  can assume
that the sum of the coefficients on  each lines is at least 2.  Let $x$
be  a  {\cita  Perron-Frobenius}  eigenvector  with  positive  entries
associated  to the  {\cita Perron-Frobenius}  eigenvalue  $\rho(V)$ of
$V$.  Let  $i_0$ be such that $\displaystyle  \min_{j=1,\dots,d} x_j =
x_{i_0} >0$. Then
$$  \rho(V)  x_{i_0}  =   \sum_{j=1}^d  v_{i_0  j}  x_j  \geq  x_{i_0}
\sum_{j=1}^d v_{i_0j}.
$$ Since $x_{i_0}>0$, we see that there exists $i_0$ such that
$$ \rho(V) \geq \sum_{j=1}^d v_{i_0j}\geq 2.
$$
\end{proof}

\begin{Proposition}\label{prop:fix:sing}
Let $\phi$  be a pseudo-{\cita Anosov} homeomorphism  affine on a translation
surface  $(X,\omega) \in  \mathcal{H}^{hyp}(2g-2)$ or  $(X,\omega) \in
\mathcal{H}^{hyp}(g-1,g-1)$.  If $\phi$  fixes a horizontal separatrix
emanating from a singular point of $\omega$  then the dilatation of
$\phi$ is bounded from below by $2$.
\end{Proposition}

\begin{proof}
By {\cita Veech}'s Theorem, $\phi$ is obtained by taking a closed loop
in a {\cita  Rauzy} diagram.  By a Lemma of  {\cita Rauzy} (see~\cite{Rauzy}), any {\cita
Rauzy} diagram contains a vertex of the kind:
$$\nu=\left(\begin{array}{ccc}   \alpha&\ldots&\beta  \\   \beta  &\ldots&
\alpha
\end{array}\right).$$
Since the  underlying flat surface  is in the  hyperelliptic connected
component,  it is  easy  to see  that  $\nu=\tau_n$ and therefore this {\cita  Rauzy} diagram  is
necessarily $\mathcal{D}^{hyp}$ (one can  also use the main theorem of
\cite{Boissy:rc}).  Has  we  have  seen,  the  labeled  and  reduced
diagrams coincide. Hence, for any closed path in the reduced diagram
whose associated matrix is  primitive, the corresponding lift contains
all the letters as winner. Hence by Proposition~\ref{prop:win:los} the associated dilatation is at least
$2$.
\end{proof}

\subsection{Case when the fixed separatrix is adjacent to a regular point}

One needs to prove the following proposition.

\begin{Proposition}\label{prop:fix:reg}
Let $\phi$  be a pseudo-{\cita  Anosov} homeomorphism on a  surface in
$\mathcal{H}^{hyp}(2g-2)$ or  $\mathcal{H}^{hyp}(g-1,g-1)$.  If $\phi$
fixes a regular  point, then the dilatation of  $\phi$ is bounded from
below by $2$.
\end{Proposition}

The   idea    of   the    proof   is   similar    to   the    one   of
Proposition~\ref{prop:fix:sing}:   we   must  consider   pseudo-{\cita
Anosov}  homeomorphisms  that  are   obtained  by  the  {\cita  Veech}
construction using the  diagram $\mathcal D_{r}(\pi_{n})$, were $\pi_n
= \left( \begin{smallmatrix}  0 & 2 & 3 &  \dots & n-1 & 1 &  n \\ n &
n-1 & \dots & 3 & 2  & 1 & 0 \end{smallmatrix} \right)$. In this case,
the  reduced diagram  is  different from  the  labeled diagram  (see
section~\ref{sec:description:marked}).  But we will  show that  we can
still apply Proposition~\ref{prop:win:los}.

\begin{Remark}
Our case is more subtle than the previous case. Indeed,
start with the permutation $\left( \begin{smallmatrix} 1 & 3 & 0 & 2 &
4  \\ 4  &  3 &  2  & 1  & 0  \end{smallmatrix}  \right) \in  \mathcal
D(\pi_{4})$. If we consider the path $b-t-b-t-t-t$ then the letter $3$
is never winner  nor looser.  But the corresponding  path in $\mathcal
D_{r}(\pi_{4}) $ is closed and primitive.  Nevertheless the lift based
at the permutation $\left( \begin{smallmatrix} 1 &  2 & 3 & 0 & 4 \\ 4
& 3 &  2 & 0 & 1 \end{smallmatrix} \right)$  is $t-t-t-b-t-b$, and all
the letters are looser.
\end{Remark}

We  will call the  permutations in  the diagram  $\mathcal D(\pi_{n})$
that   are  in   $\{\pi_{n},\pi_n\circ   \sigma,  \ldots,   \pi_n\circ
\sigma^{n-2}\}$  the   {\it  central  permutations}.    We  will  call
permutations in the diagram  $\mathcal D(\pi_{n})$ that corresponds to
$\mathcal A_{n}$ up to renumbering with $\sigma^k$ the {\it transition
permutations}.

We have the following technical lemmas.

\begin{Lemma}
\label{lm1:annexe}
Let $\gamma$ be a primitive closed path in $\mathcal D_{r}(\pi_{n})$.
There is  a lift $\hat  \gamma$ of $\gamma$ in  $\mathcal D(\pi_{n})$,
not necessarily closed, that starts and ends at central permutations.
\end{Lemma}

\begin{Lemma}
\label{lm2:annexe}
Let $\hat \gamma$ be a  path in $\mathcal{D}_r(\pi_n)$ that starts and
ends  at central permutations.   Then $\hat  \gamma$ contains  all the
letters $0,1,2,\dots,n$ as looser.
\end{Lemma}

We first prove the proposition assuming Lemma~\ref{lm1:annexe} and Lemma~\ref{lm2:annexe}.

\begin{proof}[Proof of Proposition~\ref{prop:fix:reg}]
Let $\phi$ be a pseudo-{\cita Anosov} homeomorphism that fixes a regular point. Then
it  also  fixes its horizontal outgoing separatrix. 
By {\cita Veech}'s theorem $\pA$ is obtain by taking a closed loop in some 
{\cita Rauzy} diagram.  Using {\cita Rauzy} Lemma stated  in the proof
of Proposition~\ref{prop:fix:sing}, we can show that the {\cita Rauzy}
diagram is  necessarily $\mathcal{D}_r(\pi_n)$  (one can also  use the
main theorem of \cite{Boissy:rc}).
Let   $\gamma$  be   the  corresponding   closed  primitive   path  in
$\mathcal{D}_r(\pi_n)$. By Lemma~\ref{lm1:annexe},  there is a lift of
$\gamma$  that starts  and  ends at  central  permutations. This  lift
contains all the letters  as loosers by Lemma~\ref{lm2:annexe}, and by
Proposition~\ref{prop:win:los},     the    {\cita    Perron-Frobenius}
eigenvalue of  the corresponding matrix,  and hence the  dilatation of $\pA$
is at least $2$.
\end{proof}

The proof  of the lemmas  is strongly related  to the geometry  of the
diagram $\mathcal{D}(\pi_n)$.  Before giving a formal  proof,
we first present a very informal proof that uses the informal description
of   the   diagram   given   in   the   beginning   of   Section~\ref{sec:description:marked}.
\medskip

\begin{enumerate}
\item  For the  first lemma  it is  enough to  prove that  a primitive
closed path in  the labeled {\cita Rauzy} diagram  must pass through a
central  permutation.  Such path  must pass  through a  principal loop
were  0 or  n  is  winner.  If  it  enters in  a  k-th secondary  loop
(attached to  the k-th vertex of  a principal loop) ,  the geometry of
the diagram imposes that it  either leaves the secondary loop from the
same vertex (and therefore joins  the (k+1)-th vertex of the principal
loop) or joins a k-th secondary  loop in another leaf and escape it at
the k-th  vertex of  another principal loop,  and therefore  joins the
(k+1)-th vertex  of this loop. Iterating this  argument, it eventually
joins a central permutation.
\item The key observation for the  second lemma is that a path joining
two central permutations will either  pass through all the vertices of
a principal loop, or will pass  from one principal loop to another one
through  a  transition  permutation  and  the  corresponding  pair  of
secondary  loops. In  any cases  all the  letters  will appear  as
loosers.
\end{enumerate}

We now give a proof of the lemmas.

\begin{proof}[Proof of Lemma~\ref{lm1:annexe}]
Let     $\gamma$     be     a    primitive    closed     path     in
$\mathcal{D}_{r}(\pi_{n})$,  and  let  $\hat  \gamma$  be  a  lift  of
$\gamma$  in  $\mathcal{D}(\pi_n)$.  The  path $\hat  \gamma$  is  not
necessarily closed,  but a power  of $\gamma$
admits  a lift  $\eta$  which is  closed  and primitive.  Furthermore,
$\eta$  consists  of  a   concatenation  of  lifts  of  $\gamma$.   By
Proposition~\ref{prop:closed} the path  $\eta$ contains all the letter
$0,\dots,n$   as   winner,   in   particular   the   letter   $0$   as
winner. \medskip

The point is  that the letter $0$ appears as winner  only on the steps
of  the form  $\mathcal{R}_{b}^{k-1}  (\pi_n\circ \sigma^i)\rightarrow
\mathcal{R}_{b}^{k}(\pi_n\circ\sigma^i)$.  One can assume  that $i=0$.
By   Lemma~\ref{subpath},   there    exists   in   $\eta$   the   step
$\mathcal{R}_b^{k}(\pi_n)                             \xrightarrow{0,*}
\mathcal{R}_b^{k+1}(\pi_n)$ or  the step $\mathcal{R}_t^{k}(\pi_n\circ
\sigma^{j})\xrightarrow{n,*}             \mathcal{R}_t^{k+1}(\pi_n\circ
\sigma^{j})$,   where  $j$   corresponds  to   the  index   such  that
$\mathcal{R}_b                               \mathcal{R}_t^k(\pi_n\circ
\sigma^j)=\mathcal{R}_t\mathcal{R}_b^k(\pi_n)$.

Since, $0$ and $n$ play a  symmetric role, we can iterate the argument
and therefore,  we must  reach $k=n$, which  corresponds to  a central
permutation. Hence $\eta$ contains a central permutation. Since $\eta$
is a concatenation of lifts  of $\gamma$, this path passes through the
reduced permutation $\pi_n$. Since $\gamma$ is closed, we can assume
that it starts  and ends at $\pi_n$. Hence, any  lift of $\gamma$ has
endpoints that are  in the preimage of $\pi_n$,  which are the central
permutations.  Lemma~\ref{lm1:annexe} is proven.
\end{proof}

\begin{proof}[Proof of Lemma~\ref{lm2:annexe}]
Let  $\hat \gamma$  be a  path in  $\mathcal D(\pi_n)$  connecting two
central permutations. For simplicity,  assume that $\hat \gamma(1)$ is
$\pi_n$.  Again, by symmetry,  one can  assume that  the first
arrow is given by the map $\mathcal R_{t}$, \emph{i.e.} $\hat \gamma(1)
\xrightarrow{n,0}     \hat     \gamma(2)=\mathcal{R}_t(\pi_n)$.     By
Lemma~\ref{subpath} there  exists $i_2,\ldots ,i_{n-1}$  such that the
path        $\hat       \gamma$        contains        the       steps
$\mathcal{R}_{\varepsilon_2}^2(\pi_n   \circ  \sigma^{i_2})\rightarrow
\mathcal{R}_{\varepsilon_2}^3(\pi_n   \circ   \sigma^{i_2})   ,\ldots,
\mathcal{R}_{\varepsilon_{n-1}}^{n-1}(\pi_n                       \circ
\sigma^{i_{n-1}})\rightarrow \mathcal{R}_{\varepsilon_{n-1}}^{n}(\pi_n
\circ \sigma^{i_n})$.

The first possibility  is that $\hat \gamma$ does  not change of leaf,
\emph{i.e.} for  all $k\in \{2,\ldots,n-1\}$, $i_k=0$.  Then the steps
previously  written form  a subpath  that  can be  rewritten as  $\hat
\gamma(1)  \xrightarrow{n,0} \hat  \gamma(j_2)  \xrightarrow{n,1} \hat
\gamma(j_3)    \xrightarrow{n,2}   \dots    \xrightarrow{n,n-1}   \hat
\gamma(1)$. Then  all the letters $0,\dots,n-1$ are  loosers. The path
$\hat \gamma$ is  closed, so $n$ cannot always be  winner. Hence it is
also looser.
\medskip

The  second  possibility  is  that  $\hat  \gamma$  changes  of  leaf,
\emph{i.e.} there exists a smallest $k$ such that $i_k > 0$. Then by
Lemma~\ref{subpath},  a  subpath  of   $\gamma$  is  obtained  by  the
following way (compare with Figure~\ref{fig:diagram:labeled}):
\begin{itemize}
\item We start from $\pi_n$.
\item We apply $k$ times the {\cita Rauzy} move $\mathcal{R}_t$ (these
are  the  moves  $\hat  \gamma(1) \xrightarrow{n,0}  \hat  \gamma(j_2)
\xrightarrow{n,1}         \dots         \xrightarrow{n,k}         \hat
\gamma(j_{k+1})=\mathcal{R}_t^k(\pi_n)$).
\item We apply $1$ time the {\cita Rauzy} move $\mathcal{R}_b$ (we reach
the  permutation  $\mathcal{R}_b  \mathcal{R}_t^k(\pi_n)=\mathcal{R}_t
\mathcal{R}_b^k(\pi_n\circ\sigma^{i_k})$).   This   is   the  move   $
\gamma(j_{k+1}) \xrightarrow{k,n} \hat \gamma(j_{k+2})$.
\item We  apply $k'=n-k$ times the {\cita  Rauzy} move $\mathcal{R}_t$
(until   we    get   the   permutation   $\mathcal{R}_{b}^k(\pi_n\circ
\sigma^{i_k})$).   These   are   the   moves   $\hat   \gamma(j_{k+2})
\xrightarrow{1,k}   \hat  \gamma(j_{k+3})   \xrightarrow{1,k+1}  \dots
\xrightarrow{1,n-1} \hat \gamma(j_{n+2})$.
\end{itemize}

We   see   that   all   the   letters   are   loosers.   This   proves
Lemma~\ref{lm2:annexe}.
\end{proof}

\appendix

\section{Examples in hyperelliptic components}
\label{appendix:examples}

In  this appendix  we  show that  the  uniform lower bound  on dilatations  of
pseudo-{\cita            Anosov}           homeomorphisms           in
Theorem~\ref{theo:main:general}  is  sharp  by  constructing  suitable
examples.      This     will      thus     give     a     proof     of
Theorem~\ref{theo:main:sharp}.

\subsection{Hyperelliptic connected component $\mathcal H^{hyp}(2g-2)$}
\begin{Proposition}
Let $g\geq 2$. There  exists  a  pseudo-{\cita  Anosov} homeomorphism  $\pA_g$  affine on  a translation 
surface  in  $\mathcal{H}^{hyp}(2g-2)$   whose  dilatation is the  {\cita
Perron} root  of the polynomial $X^{2g+1}  - 2X^{{2g-1}} -  2X^2 + 1$.
This dilatation satisfies
$$ 0 < \dil(\pA_g) - \sqrt{2} < \frac1{2^{2g-3}}.
$$

\end{Proposition}

\begin{Lemma}
The  polynomial $P_{g}=X^{2g+1}  - 2X^{{2g-1}}  - 2X^2  + 1$  admits a
unique real root $\alpha$ greater than $\sqrt{2}$.
\end{Lemma}
\begin{proof}
$P_g(\sqrt{2})=-4+1<0$, and $P_g(x)>0$ for $x$ large enough and $P_g'(x)>0$ for $x>\sqrt{2}$.
\end{proof}

\begin{Lemma}
Let $M=[m_{ij}]$ be the matrix of dimension $2g$ defined by:
\begin{itemize}
\item for all $j\in \{g+1,\ldots,2g\}$, $m_{1,j} =1$, and $m_{1,g}=2$,
\item for all $i\in \{2,\ldots,g\}$, $m_{i,i+g-1}=1$,
\item for all $j\in \{g+1,\ldots,2g-1\}$, $m_{i,i-g}=1$,
\item $m_{2g,g}=m_{2g,2g}=1$.
\item all the other elements are zero.
\end{itemize}
$$ M_g= \left(
\begin{array}{c|c|c}
0_{g\times g-1} &
\begin{array}{ccccc}
2&1&\cdots&\cdots&1\\  0&1&0&\cdots&0\\  \vdots  &\ddots  &  \ddots  &
\ddots& \vdots\\ \vdots & & \ddots & \ddots &0 \\ 0& \cdots & \cdots&0
&1
\end{array}& 
\begin{array}{c} 1\\ 0 \\ \vdots  \\ \vdots \\ 0 \end{array} \\
\hline I_{g-1\times g-1}& 0_{g-1\times g}&
\begin{array}{c} 0\\ \vdots\\\vdots  \\0 \end{array}\\
\hline 0& \begin{array}{ccccc} 1\ &0\ &\cdots\ &\cdots\ & 0
\end{array}  &1
\end{array}
\right)
$$
Then $M_g$ contains $\alpha$ and $1/\alpha$ as eigenvalues.
\end{Lemma}

\begin{proof}
We      must      compute      the      characteristic      polynomial
$\chi_g=\det(M_g-xI_{2g \times 2g})$. Since $P_g$  is reciprocal and $\alpha\neq -1$,
we  just need to  show that  $(X+1)\chi_g=P_g$ for  all $g\geq  2$. We
denote by  $L_i$ the $i-th$ line.  For each $i\in  \{1,g\}$ we replace
$L_i$  by $L_i+xL_{i+g}$. Then  we develop  the determinant  $g$ times
along the first column, to get:
\begin{eqnarray*}
\det(M_g-xI)&=&(-1)^{(g+2)(g-1)}\left(   2(1-x)   +  (-1)^{g+2}\left(
(-1)^{g+1}.1+  x(1-x)  D_{g-1} \right)  \right),\\  &=&2 (1-x)+  (-1)+
(-1)^gx(1-x)D_{g-1},
\end{eqnarray*}
were $D_{g-1}$ is the $(g-1)\times(g-1)$ determinant:
$$D_{g-1}=\left|   \begin{array}{cccccc}   1-x^2&1&1&\ldots&\ldots&1\\
1&-x^2&0&\ldots&\ldots&0\\            0&1&-x^2&0&\ldots&0           \\
\vdots&\ddots&\ddots&\ddots&\ddots&\vdots  \\ \vdots &&0&1&-x^2&  0 \\
0&\ldots&\ldots&0&1&-x^2
\end{array}\right|.
$$    Developing   $D_g$    on    the   last    column,   we    obtain
$D_g=-x^2D_{g-1}+(-1)^{g+1}$.  Then using this last expression, we see
that:
$$(X+1)\left(\chi_{g+1}-X^2\chi_{g}\right)=(X+1)(1-X-2X^2+2X^3)=(2X^4-3X^2+1)=P_{g+1}-X^2
P_g$$  Since $(1+X)\chi_2=P_2$, we  deduce that  $(1+X)\chi_g=P_g$ for
all $g\geq 2$, which proves the proposition.
\end{proof}

\begin{Lemma}
Let $\lambda$ (respectively, $\tau$) be  an eigenvector of $M$ for the
eigenvalue $\alpha$ (respectively,  $1/\alpha$). Then, up to replacing
$\lambda$  or  $\tau$  by  its opposite,  $(\lambda,\tau)$  defines  a
suspension datum for the permutations $\pi$ and $\pi'$, with:
$$ \pi=\left(\begin{smallmatrix} 1&2&\dots&\ldots&\ldots&2g-1&2g\\ 2g&
g&\ldots&1&2g-1&\ldots&g+1
\end{smallmatrix}\right)
\qquad   \textrm{and}  \qquad  \pi'=\left(   \begin{smallmatrix}  g+1&
\ldots&        \ldots&         2g-1&        1&\ldots&g&2g\\        2g&
2g-1&\ldots&\ldots&\ldots&\ldots& 2 & 1
\end{smallmatrix}\right)
$$ Moreover $(\pi,\lambda,\tau)$ and $(\pi',\lambda,\tau)$ define two 
isometric surfaces.
\end{Lemma}

\begin{proof}
We            have            $M_g\lambda=\alpha\lambda$           and
$M_g\tau=\frac{1}{\alpha}\tau$. This gives:
\begin{eqnarray}
\forall i\in \{1,\ldots,2g-2\}\backslash \{g\} \quad\lambda_i=\alpha^2
\lambda_{i+1}  \quad &\textrm{and}&  \quad  \tau_i =\frac{1}{\alpha^2}
\tau_{i+1}\label{l1}\\
\lambda_{2g-1}=\alpha\lambda_{g}       \quad       &\textrm{and}&\quad
\tau_{2g-1}=\frac{1}{\alpha}\tau_{g}\label{l2}\\
\lambda_{g}+\lambda_{2g}=\alpha\lambda_{2g}  \quad &\textrm{and}&\quad
\tau_{g}+\tau_{2g}=\frac{1}{\alpha}\tau_{2g}\label{l3}
\end{eqnarray}
Since $\alpha>1$, Equation~(\ref{l3}) implies that $\lambda_{g}$ is of
the same  sign as $\lambda_{2g}$.  The other equations imply  that all
the $\lambda_i$  are of the same  sign. Hence all  the $\lambda_i$ are
positive if we choose $\lambda_{2g}>0$.

By  a  similar  argument,  for  $\tau$,  we see  that  we  can  choose
$\tau_{2g}<0$ and  have $\tau_{i}>0$ for all  $i<2g$. Furthermore, the
label ``2g''  is the  last of the  first line of  $\pi$ (respectively,
$\pi'$)  and the  first of  the  second line  of $\pi$  (respectively,
$\pi'$). Then in order to  prove that $(\lambda,\tau)$ is a suspension
data   for   $\pi$  and   $\pi'$,   it   is   enough  to   show   that
$s=\sum_{i=1}^{2g}  \tau_{i}$  is  negative.  By  Equations~(\ref{l1})
and~(\ref{l2}), we have
$$s=\tau_{2g}+\sum_{k=0}^{2g-2}\frac{1}{\alpha^k}\tau_{g}=\frac{1}{\alpha}
\tau_{2g}+\sum_{k=1}^{2g-2}\frac{1}{\alpha^k}\tau_{g}=\ldots=\frac{1}{\alpha^{2g-1}}\tau_{2g}<0.$$

Hence    we     have    proven    that     $(\pi,\lambda,\tau)$    and
$(\pi',\lambda,\tau)$  define  translation   surfaces  by  the  {\cita
Veech}'s  construction   (see  Section~\ref{sec:suspension}).  Now  we
remark    that,   if    $P_1$   is    the   polygon    associated   to
$(\pi,\lambda,\tau)$,   and  $P_2$  is   the  polygon   associated  to
$(\pi',\lambda,\tau)$,  then   we  obtain   $P_2$  from  $P_1$   by  a
$180^\circ$ rotation.  Since these flat surfaces  are hyperelliptic we
conclude  that the  two surfaces  defined by  $(\pi,\lambda,\tau)$ and
$(\pi',\lambda,\tau)$ are isometric.
\end{proof}

\begin{Lemma}
There is a pseudo-{\cita  Anosov} homeomorphism on the surface defined
by the data $(\pi,\lambda,\tau)$ with dilatation $\alpha$.
\end{Lemma}

\begin{proof}
We  start  from  $\pi'$,  and   consider  the  path  $\gamma$  in  the
corresponding {\cita  Rauzy} graph obtained  by applying to  $\pi$ map  $\mathcal{R}_b$ $g$ 
times and  then  one time  the  map  $\mathcal
R_t$. We obtain the permutation
$$\pi''=\left(   \begin{smallmatrix}   g+1&   \ldots&  \ldots&   2g-1&
1&\ldots&g&2g\\ 2g& 1&2g-1&\ldots&\ldots&\ldots&\ldots& 2
\end{smallmatrix}\right)
$$  The sequence  of winners/loosers  is  $(1,2g), (1,g),\ldots,(1,2),
(2g,1)$. And therefore, the corresponding transition matrix is:
$$ M'_g= \left(
\begin{array}{c|c|c}
\begin{array}{ccccc}
2&1&\cdots&\cdots&1\\  0&1&0&\cdots&0\\  \vdots  &\ddots  &  \ddots  &
\ddots& \vdots\\ \vdots & & \ddots & \ddots &0 \\ 0& \cdots & \cdots&0
&1
\end{array}& 
0_{g\times g-1} &
\begin{array}{c} 1\\ 0 \\ \vdots  \\ \vdots \\ 0 \end{array} \\
\hline 0_{g-1\times g}& I_{g-1\times g-1}&
\begin{array}{c} 0\\ \vdots\\\vdots  \\0 \end{array}\\
\hline
\begin{array}{ccccc} 1\ &0\ &\cdots\ &\cdots\ & 0 \end{array} 
&0 &1
\end{array}
\right)
$$

The    permutation   $\pi''$    is   obtained    from    $\pi$   after
renumbering. Multiplying $M'$ by the corresponding permutation matrix,
one  get   precisely  the  matrix  $M_{g}$.  Therefore,   if  we  consider
$\lambda''=\frac{1}{\alpha}\lambda$, renumber  it and do  backward the
path $\gamma$,  one get $(\pi',\lambda)$. Therefore, if  we start from
$(\pi,\lambda,\tau)$,  then consider  $(\pi',\lambda,\tau)$,  and then
apply the  {\cita Rauzy} induction, and renumbering,  one will obtain
$(\pi,\frac{1}{\alpha}\lambda,\alpha \tau)$. Each  time, we obtain new
parameter for  the same  flat surface, hence  $(\pi,\lambda,\tau)$ and
$(\pi,\frac{1}{\alpha},\lambda,\alpha\tau)$ define the same element in
the  moduli   space.  So  the  corresponding  flat   surface  admit  a
pseuso-Anosov homeomorphism with dilatation $\alpha$.
\end{proof}

\begin{proof}[Proof of the proposition]
The existence  of the  pseudo-{\cita Anosov} homeomorphism  $\pA_g$ is
clear from the previous lemmas.  We have $\alpha=\dil(\pA_g)$, and:
$$   \alpha^{2g+1}   -   2\alpha^{{2g-1}}    -   2\alpha^2   +   1   =
\alpha^2\left(\alpha^{2g-1} - 2\alpha^{2g-3} - 2\right) + 1 = 0,
$$ thus
$$  \alpha^{2g-1}   -  2\alpha^{2g-3}  -  2   =  \alpha^{2g-3}  \left(
\alpha^{2} - 2 \right) - 2 < 0.
$$  Since $\alpha  > \sqrt{2}$  we  obtain $0  < \alpha  - \sqrt{2}  <
\frac1{2^{2g-3}}$.

\end{proof}

\subsection{Hyperelliptic connected component $\mathcal H^{hyp}(g-1,g-1)$}
\begin{Proposition}
Let $g\geq 2$. There  exists a pseudo-{\cita Anosov}  homeomorphism  on  a translation surface  in
$\mathcal H^{hyp}(g-1,g-1)$  with dilatation  the  {\cita Perron}  root of  the
polynomial:
$$
\begin{array}{ll}
X^{2g+2} -  2 X^{2g}  - 2 X^{g+1}  - 2  X^2 + 1,  & \textrm{if  $g$ is
even,} \\  X^{2g+2} - 2 X^{2g} -  4 X^{g+2} + 4  X^{g} + 2 X^2  - 1, &
\textrm{if $g$ is odd.}
\end{array}
$$ The dilatation satisfies
$$
0 <  \dil(\pA_g) - \sqrt{2} <  \frac{4}{\sqrt{2}^{g}}
$$
\end{Proposition}

\begin{proof}
The idea  of the proof  is very similar  to the previous one.  We just present
here the corresponding paths in the {\cita Rauzy} diagram.
\begin{itemize}
\item If $g$ is even, we start from the permutation $\pi=\left(
\begin{smallmatrix}
1&\ldots &g+1&g+2&\ldots&2g&2g+1\\ 2g+1&g+1&\ldots&1&2g&\ldots&g+2
\end{smallmatrix}\right)$. As previously, we consider the ``reverse'' permutation 
$\pi'=\left(
\begin{smallmatrix}
g+2&\ldots&2g&1&\ldots&g+1&2g+1\\ 2g+1&2g&\ldots&g+2&g+1&\ldots&1
\end{smallmatrix}\right)$. 
Then the path in the {\cita Rauzy} diagram starting from $\pi'$ and defined by $(g+1)$ times $b$, two times $t$. One gets the permutation $\pi''=\left(
\begin{smallmatrix}
g+2&\ldots&2g&1&&2&3\ldots&g+1&2g+1\\
2g+1&1&2&2g&\ldots&g+2&g+1&\ldots&3
\end{smallmatrix}\right)$, which is a renumbering of $\pi$, and the associated matrix is:
$$ N_g= \left(
\begin{array}{c|c|c}
0_{g+1\times g-1} &
\begin{array}{ccccc}
2&2&1&\cdots&1\\  0&1&0&\cdots&0\\ \vdots &\ddots  & \ddots  & \ddots&
\vdots\\ \vdots & & \ddots & \ddots &0 \\ 0& \cdots & \cdots&0 &1
\end{array}& 
\begin{array}{c} 1\\ 0 \\ \vdots  \\ \vdots \\ 0 \end{array} \\
\hline I_{g-1\times g-1}& 0_{g-1\times g+1}&
\begin{array}{c} 0\\ \vdots\\\vdots  \\0 \end{array}\\
\hline 0& \begin{array}{cccccc} 1&1&0&\cdots&\cdots& 0 \end{array} &1
\end{array}
\right)
$$

One can check that $(X+1)\chi_{N_g}(X)=X^{2g+2} - 2 X^{2g} - 2 X^{g+1}
- 2  X^2 +  1$, and  that,  for $\alpha$  the {\cita  Perron} root  of
$\chi_{N_g}(X)$, the  eigenvectors of $N_g$  corresponding to $\alpha$
and  $\frac{1}{\alpha}$  define  a  suspension  data  for  $\pi$  (and
$\pi'$).  Then  the   corresponding  translation  surface  admits  the
required pseudo-{\cita Anosov} homeomorphism by construction. \medskip

Also, we have
$$ \alpha^{2g}(\alpha^2-2)< 2\alpha^{g+1} +2\alpha^2
$$ Thus,
$$              \alpha-\sqrt{2}<             \frac{2}{\alpha+\sqrt{2}}
\left(\frac{1}{\alpha^{g-1}}+\frac{1}{\alpha^{2g-2}}\right)<\frac{2}{\sqrt{2}^{g}},
$$ since $\alpha>\sqrt{2}$.
\item
If $g$  is odd, we consider $\pi$  as previously and we  take the path
defined  by $g$  times $b$,  then  $t-b-t-t$. We  obtain the  required
pseudo-{\cita  Anosov}   with  the  same   construction  we  presented have
above. Also, we have:
$$                                                     \alpha-\sqrt{2}=
\frac{1}{\alpha+\sqrt{2}}\frac{4\alpha^g(\alpha^2-1)+1-2\alpha^2}{\alpha^{2g}}<\frac{4}{\alpha^g}<\frac{4}{\sqrt{2}^g}$$
\end{itemize}

\end{proof}

\section{Examples in non hyperelliptic components}
\label{appendix:examples:2}

In      this     appendix      we     motivate      assumptions     of
Theorem~\ref{theo:main:general}.  We   show  that  if   we  relax  the
condition   on  the  hyperelliptic   involution,  one   can  construct
pseudo-{\cita   Anosov}  homeomorphisms   in  the   non  hyperelliptic
connected  components whose  dilatations tend  to $1$  when  the genus
tends   to    infinity.    This   will   thus   give    a   proof   of
Theorem~\ref{theo:main:odd}. \medskip

We will  construct a sequence of  pseudo-{\cita Anosov} homeomorphisms
$\varphi_{g}$  on the flat  sphere i.e.  $\varphi_{g}$ is  affine with
respect  to  a  quadratic  differential  $q_{g}$ on  the  sphere.  The
construction involves the  {\cita Rauzy-Veech} induction for quadratic
differentials~\cite{BL}, using \emph{generalized permutations}.

\begin{Proposition}
\label{prop:construc:append}
For each $g\geq 3$, there exists a pseudo-{\cita Anosov} homeomorphism
$\varphi_{g}$ on  the flat  sphere with $2g$  poles and two  zeroes of
order   $g-2$,   fixing   a   pole,   and   having   for   dilatations
$\dil(\varphi_{g})$ the {\cita Perron} root of the polynomial
$$ P = X^{2g} - X^{2g-1} - 4X^g - X + 1.
$$
\end{Proposition}

\begin{proof}[Proof of the proposition]
We consider the {\cita Rauzy}  diagram corresponding to the stratum of
quadratic  differentials on  the sphere  having the  given singularity
data. It  is sufficient to give a  closed path in this  diagram and to
check  that   the  renormalisation  matrix  is   irreducible  and  has
$\dil_{g}$ for eigenvalue, where $\dil_{g}$ is the {\cita Perron} root
of $P$. \medskip

Let  us  consider the  following  generalized  permutation (on  $2g+1$
letters)
$$     \left(\begin{array}{cccccccc}     1&2&2&3&3&\dots&g+1&g+1    \\
g+2&g+3&g+3&\dots&2g+1&2g+1&g+2&1 \\
\end{array}\right).
$$ For instance for $g=3$ this gives
$$ \left(\begin{array}{cccccccc} 1&2&2&3&3&4&4 \\ 5&6&6&7&7&5&1
\end{array}\right).
$$ The {\cita Rauzy} path we will consider is $t-b-t-b-t-b$.  A simple
calculation shows  that the renormalization  matrix $\widehat{V}$ (for
the labeled permutation) and permutation matrix $P$ are, respectively
$$
\begin{array}{lll}
\widehat{V} = \left(
\begin{array}{ccccccc}
2 & &&&&& 1 \\ & 1 && \\ & & \ddots & \\ 2 & & & 1 & 1 &&2\\ & & & 1 &
2 \\ & & & & & \ddots \\ 1 & & & 1 & & & 1 \\
\end{array}
\right)
& \textrm{ and } &
P = \left(
\begin{array}{c|c}
\begin{array}{c} 0 \\ \vdots \\ 0 \end{array}
& I_{2g\times2g} \\ \hline 1 & \begin{array}{ccc} 0 & \cdots & 0
\end{array}
\end{array}
\right),
\end{array}
$$ where the ``2'' in the diagonal of the matrix $\widehat{V}$ appears
at  the position  $g+2$. All  the other  entries in  $\widehat{V}$ are
zeroes.   Thus   the  renormalization   matrix   (for  the   reduced
permutations)  is $V  = \widehat{V}\cdot  P$.  One can  show that  the
matrix $V$ is irreducible and its characteristic polynomial is
$$ (X-1)(X^{2g} - X^{2g-1} - 4X^g - X + 1) = (X-1)P(X).
$$ Then as  in the previous section, we see  that some eigenvectors of
$V$ define suspension  data, and the corresponding flat  surface has a
pseudo-{\cita  Anosov} homeomorphism whose dilatation is the  {\cita Perron}  root of $P$.
\end{proof}

Finally  we   have  the  following  corollary,   which  justifies  the
assumptions  of Theorem~\ref{theo:main:general}  on  the hyperelliptic
involution

\begin{Corollary}
For $g\geq 3$ odd,  there exists an affine pseudo-{\cita Anosov} homeomorphism
$\pA_{g}$  on a  translation surface  $(X,\omega)$ with  the following
properties:
\begin{enumerate}
\item $(X,\omega) \in \HA^{odd}(g-1,g-1)$,
\item $(X,\omega)$ is hyperelliptic,
\item  $\pA_{g}$ fixes  a separatrix  on $(X,\omega)$  (issued  from a
{\cita Weierstrass} point),
\item the  dilatation of $\pA_{g}$ is  the {\cita Perron}  root of the
polynomial
$$ X^{2g} - X^{2g-1} - 4X^g - X + 1.
$$
\end{enumerate}
In particular
$$ \lim\limits_{g \to +\infty} \dil(\pA_{g}) = 1.
$$
\end{Corollary}

\begin{Remark}
Obviously the hyperelliptic involution on  $X$ fixes the two zeroes of
$\omega$, and $(X,\omega)$ is not in the hyperelliptic connected component of the corresponding stratum.
\end{Remark}

\begin{proof}[Proof of the corollary]
For       $g$       odd,       the       examples       given       by
Proposition~\ref{prop:construc:append}  lift to  pseudo-{\cita Anosov} homeomorphisms
$\pA_{g}$ on  the orientating cover $(X,\omega)  \in \HA(g-1,g-1)$. By
construction  the surface  is hyperelliptic  and since  the  poles are
ramification points,  there is  one lift that  fixe a  regular ({\cita
Weierstrass})   point   and    the   separatrix   issued   from   that
point. \medskip

To   identify   the   connected   component,  we   use   the   formula
in~\cite{Lanneau:spin}.   Since   $(X,\omega)   \rightarrow   (\mathbb
P^{1},q)$ is  the orientating cover, the spin  structure determined by
$(X,\omega)$ is
$$ \left[\cfrac{|n_{+1}-n_{-1}|}{4}\right] \qquad \mod 2,
$$ where  $n_{+1}$ is  the number of  singularities of $q$  of degrees
$1\mod  4$ and  $n_{-1}$  is the  number  of singularities  of $q$  of
degrees $-1\mod 4$.
\begin{enumerate}
\item If $g=1 \mod 4$: then  $g-2 = -1$ modulo $4$, thus $n_{+1}=0$
and  $n_{-1}=2g+2$.   Hence  the  parity  of  the  spin  structure  on
$(X,\omega)$ is (with $g=1+4k$)
$$ \left[\cfrac{|n_{+1}-n_{-1}|}{4}\right]  = \frac{2g+2}{4} =  1+2k =
1\mod 2.
$$
\item If $g=-1 \mod 4$: then  $g-2 = 1$ modulo $4$, thus $n_{+1}=2$
and  $n_{-1}=2g$.    Hence  the  parity  of  the   spin  structure  on
$(X,\omega)$ is (with $g=-1+4k$)
$$ \left[\cfrac{|n_{+1}-n_{-1}|}{4}\right] =  \frac{2g-2}{4} = -1+2k =
1\mod 2.
$$
\end{enumerate}

\end{proof}


\end{document}